\documentclass{article}[10pt]
\usepackage{amsmath,amsthm,amssymb,graphicx}

 \newtheorem{theorem}{Theorem}[section]
 \newtheorem{definition}{Definition}[section]

 \newtheorem{proposition}[definition]{Proposition}
 \newtheorem{lemma}[definition]{Lemma}
 \newtheorem{corollary}[definition]{Corollary}

 \newtheorem{example}[definition]{Example}

\numberwithin{equation}{section}

\DeclareMathOperator{\Aut}{Aut}
\DeclareMathOperator{\Sym}{Sym}
\DeclareMathOperator{\diam}{diam}
\DeclareMathOperator{\Conv}{Conv}
\DeclareMathOperator{\Fix}{Fix}
\DeclareMathOperator{\Orb}{Orb}

\begin{document}

\title{Modified Hanoi Towers Groups and Limit Spaces}
\author{Shotaro Makisumi, Grace Stadnyk, Benjamin Steinhurst\thanks{Authors supported in part by the National Science Foundation through grant DMS-0505622.}}

\maketitle

\begin{center}{\small
Contacts:\\
makisumi@princeton.edu\\
Shotaro Makisumi\\
Department of Mathematics\\
Princeton University\\
Princeton, NJ 08544 USA
\vskip 3mm
gstadnyk@hamilton.edu\\
Grace Stadnyk\\
Department of Mathematics\\
Hamilton College\\
Clinton, NY 13323 USA
\vskip 3mm
steinhurst@math.uconn.edu\\
Benjamin Steinhurst\footnote{Corresponding Author}\\
Department of Mathematics\\
University of Connecticut\\
Storrs, CT 06269 USA
}\end{center}

\begin{abstract} We introduce the $k$-peg Hanoi automorphisms and Hanoi self-similar groups, a generalization of the Hanoi Towers groups, and give conditions for them to be contractive. We analyze the limit spaces of a particular family of contracting Hanoi groups, $H_c^{(k)}$, and show that these are the unique maximal contracting Hanoi groups under a suitable symmetry condition. Finally, we provide partial results on the contraction of Hanoi groups with weaker symmetry.
\end{abstract}

\section{Introduction}
The Hanoi Towers game consists of $n$ graduated disks on $k$ pegs. The object of the game is to move all $n$ disks from one peg to another by moving one disk at a time so that no disk is on top of a smaller disk at any step. For the 3-peg case, it is well known that the optimal solution is given recursively; for this and other known results about the game, see \cite{BH99}.

This paper examines the Hanoi Towers game from two perspectives. The first is the theory of self-similar groups---groups with a \emph{self-similar action} by automorphism on regular rooted trees---a branch of geometric group theory developed in the last few decades. Any \emph{contracting} self-similar group has an associated \emph{limit space} with fractal-like properties that capture the group's self-similarity. The standard reference on self-similar groups is \cite{Nek05}. The other is analysis on fractals, which has developed analytic structures such as measures, metrics, and Laplacians on \emph{postcritically finite} (p.c.f. for short) fractals (see \cite{Kig01, Str06}) as well as the wider classes of finitely ramified and even infinitely ramified fractals, for example in \cite{Betal08,kajino,physics}. A p.c.f. fractal obtained as a limit space is thus equipped with both algebraic and analytic structures, an interplay that has driven considerable recent developments (see \cite{NT08} and references therein).

The recursive nature of the Hanoi Towers game allows us to model it using self-similar groups. In \cite{GS06} and \cite{GS08}, Grigorchuk and \v{S}un\'{i}k introduced the \emph{Hanoi Towers groups}---self-similar groups whose generators correspond to the game's legal moves---and derived some analytic results on their limit spaces. This paper aims to extend this development to modifications of the Hanoi Towers game.

After a brief introduction to self-similar groups in Section \ref{Background}, we define the $k$-peg Hanoi Towers groups, $H^{(k)}$, in Section \ref{HanoiTowersGroups} as the self-similar group generated by automorphisms $\{a_{ij}$, $0 \le i < j \le k-1\}$ each of which corresponds to a legal move between pegs $i$ and $j$ under the identification of the rooted tree with the legal states of the $k$-peg game. Though $H^{(3)}$, corresponding to the standard 3-peg game, are contracting, we find that $H^{(k)}$ for $k>3$, and in fact any interesting group generated by these automorphisms, is not contracting. We will define what we mean by ``interesting'' after Lemma \ref{noncontractioncondition}. As these non-contracting groups have no known association to fractal or self-similar limit spaces, we are led to introduce a larger class of automorphisms in Section \ref{HanoiAutomorphismsAndGroups} called the \emph{Hanoi automorphisms} from which to take generators of possibly contracting self-similar groups. We define \emph{Hanoi groups} as self-similar groups with subsets of these automorphisms as generators. The rest of the section develops sufficient conditions for these groups to be contracting; in particular, Theorem \ref{hanoicontraction2} reduces determining whether or not a Hanoi group is contracting to a finite calculation.

Section \ref{Hck} studies a particular family $H_c^{(k)}$ of contracting Hanoi groups and their limit spaces, denoted by $\mathcal{J}^{(k)}$. We show that there exist compact sets $K^{(k)}$ in $\mathbb{R}^{k+1}$ defined by an iterated function system that are homeomorphic to $\mathcal{J}^{(k)}$, and that this self-similar structure is p.c.f. Following the standard theory for p.c.f. self-similar sets \cite{Kig01}, we solve the renormalization problem for $\mathcal{J}^{(k)}$ and equip $K^{(k)}$ with a self-similar energy and effective resistance metric. From this, we calculate the Hausdorff and spectral dimensions of $K^{(k)}$.

To reflect the symmetry of pegs in the original Hanoi Towers game, we introduce a symmetry condition on the generating set of the Hanoi groups in Section \ref{Symmetry}. Under strong enough conditions, we show that $H_c^{(k)}$ is the unique maximal contracting $k$-peg Hanoi group, all other contracting Hanoi groups are subgroups of $H_c^{(k)}$. Finally, we present partial results on contracting Hanoi groups that arise when these symmetry conditions are relaxed.

The Appendix deals with Hanoi Networks, HN3 and HN4, introduced in \cite{BGG08}, which are partially inspired by the Hanoi Towers game and also possess self-similar qualities. Though they are primarily of interest in the physics literature for their small world properties which give rise to anomalous diffusion, there are some connections to the automata which arise from a different construction.

\section*{Acknowledgements}
The authors thank Alexander Teplyaev, Volodymyr Nekrashevych, Luke Rogers, Matt Begu\'{e}, Levi DeValve, and David Miller.

\section{Self-Similar Groups}\label{Background}
In this section we review the necessary background on self-similar groups. For more details, see \cite{Nek05} and \cite{NT08}.

Let $\mathsf{X}$ be a finite set, called the alphabet. Write $\mathsf{X}^n$ for the set of words of length $n$ over $\mathsf{X}$, $w = x_n \ldots x_1$, where $x_i \in \mathsf{X}$. The length of $w$ is denoted $|w|$. The union of all finite words, including the empty word $\varnothing$, is denoted $\mathsf{X}^* = \bigcup_{n = 0}^{\infty} X^n$.

The free monoid, $\mathsf{X}^*$, has a rooted tree structure with root $\varnothing$ and an edge between $w$ and $xw$ for every word $w$ and every $x \in \mathsf{X}$. An automorphism of $\mathsf{X}^*$ is a permutation that fixes $\varnothing$ and preserves adjacency; we write $\Aut \mathsf{X}^*$ for the group of these automorphisms under composition. Any automorphism acts as a permutation of the vertices of $\mathsf{X}^n$ for each $n \ge 0$.

Let $g \in \Aut \mathsf{X}^*$. Identifying $\mathsf{X}^1$ with $\mathsf{X}$, the restriction of $g$ to $\mathsf{X}^1$ becomes a permutation of $\mathsf{X}$. This is called the \emph{root permutation} of $g$, denoted $\sigma_g$. For each $x \in \mathsf{X}$, $g$ is an adjacency-preserving bijection from $x\mathsf{X}^*$ to $\sigma_g(x)\mathsf{X}^*$. Identifying both subtrees with $\mathsf{X}^*$, $g$ restricted to $x\mathsf{X}^*$ becomes another automorphism of $\mathsf{X}^*$, called the \emph{section of $g$ at $x$} and denoted $g|_{x}$. That is, we have
\begin{eqnarray} \label{sectiondefinition}
g(xw) = \sigma_g(x)g|_x(w)
\end{eqnarray}
for all $x \in \mathsf{X}$ and $w \in \mathsf{X}^*$. For $w = x_n \ldots x_1 \in \mathsf{X}^n$, $n \ge 2$, we define the section at a word inductively by $g|_w = (g|_{x_n})|_{x_{n-1}\ldots x_1}$.

The action of $g$ on $\mathsf{X}^1$ is determined by $\sigma_g$. Thereafter, the action on $\mathsf{X}^{n+1}$ is given by the action on $\mathsf{X}^n$ by (\ref{sectiondefinition}). The root permutation and the set of length-one sections uniquely determine an automorphism. Fixing $\mathsf{X} = \{0, \ldots, k-1 \}$, we use ``wreath recursion'' notation to express this dependence:
\[
g = \sigma_g (g|_0, g|_1, \ldots, g|_{k-1}).
\]
This allows us to compute compositions. Given $g$, $h \in \Aut \mathsf{X}^*$,
\[
gh(xw) = g(\sigma_h(x) h|_{x}(w)) = \sigma_g\sigma_h(x) g|_{\sigma_h(x)}h|_{x}(w)
\]
for all $x \in \mathsf{X}$ and $w \in \mathsf{X}^*$. Therefore,
\[
gh = \sigma_g\sigma_h (g|_{\sigma_h(1)}h|_{1}, \ldots, g|_{\sigma_h(k-1)}h|_{k-1}).
\]

We obtain the following lemma by applying $(gh)|_j = g|_{\sigma_h(j)}h|_j$ repeatedly to \\ $( \cdots ((a_m a_{m-1})a_{m-2}) \cdots a_2)a_1$.

\begin{lemma} \label{sectionformula}
Let $a_i \in \Aut(\mathsf{X}^{-\omega})$ with root permutation $\sigma_i$. Then for any $j \in \mathsf{X}$,
\[
(a_m \cdots a_2a_1)|_j = a_m|_{j_m} \cdots a_2|_{j_2} a_1|_{j_1},
\]
where $j_1 = j$ and $j_l = \sigma_{l - 1} \cdots \sigma_2 \sigma_1(j)$ for $l \ge 2$.
\end{lemma}

\begin{definition}
A \emph{self-similar group} $(G, \mathsf{X})$ is a group $G$ together with a faithful action of $G$ by automorphisms on $\mathsf{X}^*$ such that, for every $g \in G$ and every $x \in \mathsf{X}$, there exists some $h \in G$ such that
\[
g(xw) = g(x)h(w)
\]
for every $w \in \mathsf{X}^*$. We identify $G$ with its image in $\Aut \mathsf{X}^*$. As with tree automorphisms, we define the \emph{root permutation of $g$} and its \emph{section at $x$} by $\sigma_g = g|_{\mathsf{X}}$ and $g|_x = h$. Where there is no risk of confusion, we write $G$ for the self-similar group.
\end{definition}

The following property plays a central role in the theory of self-similar groups.

\begin{definition}
A self-similar group $(G, \mathsf{X})$ is \emph{contracting} if there exists a finite set $\mathcal{N} \subset G$ such that for every $g \in G$ there exists $k \in \mathbb{N}$ such that $g|_{v} \in \mathcal{N}$ for all $v \in \mathsf{X}^*$ with $|v| \ge k$. The smallest such $\mathcal{N}$ is called the \emph{nucleus} of $G$.
\end{definition}

Note that if $(G, \mathsf{X})$ is contracting with nucleus $\mathcal{N}$, then in particular $\mathcal{N}$ contains any $g \in G$ where $g|_i = g$ for some $i \in \mathsf{X}$.

\begin{definition} An \emph{automaton} $\mathsf{A}$ over $\mathsf{X}$ is given by
\begin{enumerate}
\item the \emph{set of states}, also denoted $\mathsf{A}$
\item a map $\tau : \mathsf{A} \times \mathsf{X} \rightarrow \mathsf{X} \times \mathsf{A}$.
\end{enumerate}
An automaton is \emph{finite} if its set of states is finite.

An automaton $\mathsf{A}$ is represented by its \emph{Moore diagram}, a digraph with vertex set $\mathsf{A}$. For each $(g, x) \in \mathsf{A} \times \mathsf{X}$, there is an edge from $g$ to $h$ labeled $(x, y)$, where $h$ and $y$ are defined by $\tau(g, x) = (y, h)$.
\end{definition}

The nucleus $\mathcal{N}$ of a contracting group can be viewed as a finite automaton with states $\mathcal{N}$ and $\tau(g, x) = (\sigma_g(x), g|_x)$. In the Moore diagram, the edge from $g$ to $g|_x$ is labeled $(x, g(x))$.

A set $S \subset \Aut \mathsf{X}^*$ of automorphisms is said to be \emph{state-closed} if every section of $g \in S$ is also in $S$. Consider a group $G$ generated by a state-closed set,$S$, of automorphisms with finite order. Then every element of $G$ can be written as a product of generators, without using inverses. Then by Lemma \ref{sectionformula}, $G$ is also state-closed and hence is a self-similar group.

\begin{example} \label{H3example}
The \emph{3-peg Hanoi Towers group} $(H^{(3)}, \{0, 1, 2\})$ is the self-similar group generated by
\[
S = \{1\} \cup \left\{
\begin{array}{c}
a_{01} = (0\ 1)(1, 1, a_{01}) \\
a_{02} = (0\ 2)(1, a_{02}, 1) \\
a_{12} = (1\ 2)(a_{12}, 1, 1)
\end{array}
\right\}.
\]

\begin{figure}[htbp]
\begin{center}
\includegraphics[width=60mm]{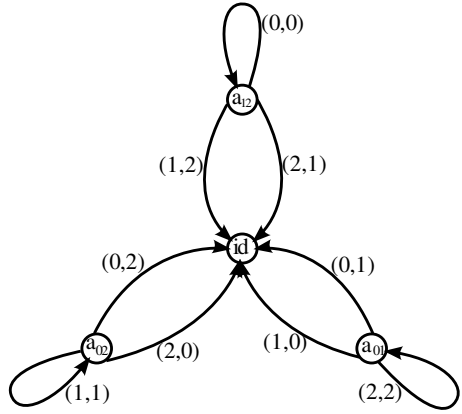}
\end{center}
\caption{Moore diagram for the nucleus of $H^{(3)}$}
\label{H3moore}
\end{figure}

Here and elsewhere, we write 1 for the identity automorphism. Note that $a_{01}^2 = a_{02}^2 = a_{12}^2$ and that, since 1 is included, $S$ is state-closed. We will later show that $H^{(3)}$ is contracting with nucleus $S$. Figure 1 shows the Moore diagram of $S$. The three self-loops of the identity automorphism, at the center, are not shown.
\end{example}

Contracting self-similar groups have an associated topological space, called their limit space, this is often also self-similar. It is this space that connects the theory of self-similar groups to fractal analysis.

\begin{definition}
Let $(G, \mathsf{X})$ be a contracting self-similar group. We write $\mathsf{X}^{-\omega}$ for the set of left-infinite sequences $\ldots x_2x_1$, $x_i \in \mathsf{X}$, and give it the product topology, where each copy of $\mathsf{X}$ has the discrete topology. Two left-infinite sequences $\ldots x_2x_1$, $\ldots y_2y_1 \in \mathsf{X}^{-\omega}$ are said to be \emph{asymptotically equivalent} if there exists a sequence $g_k \in G$, taking only finitely many different values, such that $g_k(x_k \ldots x_2x_1) = y_k \ldots y_2y_1$ for all $k \ge 1$. The quotient of $\mathsf{X}^{-\omega}$ by the asymptotic equivalence relation is called the \emph{limit space} of $G$ and is denoted $\mathcal{J}_G$.
\end{definition}

The next proposition allows us to read off asymptotically equivalent sequences from the Moore diagram of the nucleus.

\begin{proposition} (part of Theorem 3.6.3 in \cite{Nek05}) \label{nucleusasymtotic}
Two sequences $\ldots x_2x_1$, and $\ldots y_2y_1$ are asymptotically equivalent if and only if there exists a left-infinite path $\ldots e_2e_1$ in the Moore diagram of the nucleus such that the edge $e_i$ is labeled by $(x_i, y_i)$.
\end{proposition}

For example, from Figure \ref{H3moore}, we see that pairs of asymptotically equivalent sequences for $H^{(3)}$ have the form $(\ldots llli, \ldots lllj)$, where $i, j, l \in \{0,1,2 \}$ are all distinct. We can also associate a sequence of finite graphs to any self-similar group with a fixed generated set.

\begin{definition}
Let $(G, \mathsf{X})$ be a self-similar group generated by $S$. The \emph{$n$-th level Schreier graph} of $G$ with respect to $S$, denoted $\Gamma_n$, is the graph with vertices identified with $\mathsf{X}^n$ and whose vertices $w, v$ are connected if $s(w) = v$ for some automorphism $s \in S \cup S^{-1}$.
\end{definition}

When $G$ is contracting, it is known that $\{ \Gamma_n \}$ limits to $\mathcal{J}_G$ in an appropriate sense (see Section 3.6.3 in \cite{Nek05}).

\section{Hanoi Towers Groups}\label{HanoiTowersGroups}
The Hanoi Towers groups, $\{ H^{(k)}\ |\ k \ge 3 \}$, are self-similar groups introduced in \cite{GS06} to model the $k$-peg Hanoi Towers game. Fix the alphabet $\mathsf{X}_k = \{ 0, 1, \ldots, k-1 \}$, and identify it with the $k$ pegs. Consider a game with $n$-disks labeled 1 through $n$ from largest to smallest. The word $x_n x_{n-1}\ldots x_1$ uniquely determines a legal $n$-disk configuration where disk $i$ is on peg $x_i$. For each $n \ge 1$, $\mathsf{X}_k^n$ is identified with the legal states of the $n$-disk game.

Recall $a_{01} = (0\ 1)(1, 1, a_{01})$ from $(H^{(3)}, \mathsf{X}_3)$ in Example \ref{H3example}. We have $$a_{01}(x_n x_{n-1} \ldots x_1) = a_{01}(x_n)a_{01}|_{x_n}(x_{n-1} \ldots x_1).$$ If $x_n$ = 0 or 1, $a_{01}$ replaces it with the other letter, and the remaining word is unchanged; if $x_n = 2$, then $x_n$ remains unchanged, and $a_{01}$ acts on $x_{n-1} \ldots x_1$. Thus $a_{01}$ looks for the leftmost occurrence of either 0 or 1 and replaces it according to its root permutation. In the context of the game, $a_{01}$ ignores all disks on peg 2 and moves the smallest disk on either peg 0 or peg 1 to the other. This is the legal move between pegs 0 and 1.

In general, for $k \ge 3$, we define $a_{ij}$, $0 \leq i < j < k$, as the automorphism of $\mathsf{X}_k^*$ with root permutation $(i\ j)$ and
\[
a_{ij}|_{l}=\left\{
\begin{array}{cl}
1 & \mbox{if } l = i, j \\ 
a_{ij} & \mbox{otherwise.}
\end{array}
\right.
\]
As with $a_{01}$, $a_{ij}$ corresponds to the legal move between pegs $i$ and $j$. This motivates the following definition from \cite{GS06}.

\begin{definition}
The \emph{$k$-peg Hanoi Towers Group} $(H^{(k)}, \mathsf{X}_k)$, $k \ge 3$, is the self-similar group generated by $\{ a_{ij}\ |\ 0 \le i < j < k \}$.
\end{definition}

\begin{definition}\label{word length}
Given a group $G$ with a fixed generating set, $S$, a \emph{representation} of $g \in G$ is any expression $g = s_n \dots s_2 s_1$, where $s_i \in S \cup S^{-1}$. Then $n$ is called the \emph{length} of the representation. A \emph{minimal representation} is a representation of shortest length. The \emph{length of $g$}, denoted $l(g)$, is  0 if $g$ is the identity element, and otherwise is the length of a minimal representation.
\end{definition}
 
\begin{theorem}[Theorem 2.1 in \cite{GS06}] The 3-peg Hanoi Towers group, $(H^{(3)}, \mathsf{X}_3),$ is contracting.
\end{theorem}

\begin{proof}
Let $\mathcal{N} = \{1, a_{01}, a_{02}, a_{12}\} \subset H^{(3)}$. We induct on $l(g)$ to show that $l(g|_i) < l(g)$ for any $g \in H^{(3)}$ with $l(g) \ge 2$ and any $i \in \mathsf{X}_3$.

First consider the base case $l(g) = 2$. Since each $a_{ij}$ has order 2, it suffices without loss of generality to consider $a_{01}a_{12} = (1\ 2\ 0)(a_{12}, a_{01}, 1)$ all of whose sections have length less than two. Now consider $g \in H^{(3)}$ with $l(g) = n + 1$, $n \ge 2$. Without loss of generality, $g = g'a_{01}$ for some $g' \in G$ with $l(g') = n$. Then
\[
g = g'a_{01} = \sigma_{g'}(0\ 1)(g'|_1, g'|_0, g'|_2a_{01}).
\]
Assuming the claim for $n$, we have $l(g'|_1)<n$ and $l(g'|_0) < n$ and $l(g'|_2a_{01}) < n + 1$, which proves the $n + 1$ case.

Since $\mathcal{N}$ is state-closed, this shows that for any $g \in G$ with $l(g) = n$, we have $g|_v \in \mathcal{N}$ for all $v \in \mathsf{X}^*$ with $|v| \ge n - 1$. Thus $H^{(3)}$ is contracting.
\end{proof}

It happens that $H^{(k)}$ is not contracting for any $k > 3$. In fact, we can make a stronger statement.

\begin{proposition}\label{HTnoncontracting}
Any self-similar group $(G, \mathsf{X}_k)$, $k > 3$, that contains the automorphisms $a_{ij}$ and $a_{jl}$ for $i$, $j$, and $l$ all distinct is not contracting.
\end{proposition}

We first give a lemma that provides a sufficient condition for a group to not be contracting.

\begin{lemma}\label{noncontractioncondition}
Let $(G, \mathsf{X})$ be a self-similar group. Suppose some $g \in G$ satisfies the following conditions:
\begin{enumerate}
\item $\sigma_g$ is nontrivial with order $n>1$
\item There exists $i \in \mathsf{X}$ fixed by $\sigma_g$ such that $g|_i = g$.
\item There exists $j \in \mathsf{X}$ such that $g^n|_j = g^m$ for some integer $m$, $0 < |m| < n$.
\end{enumerate}
Then $G$ is not contracting.
\end{lemma}

\begin{proof}
Note that if $g\in G$ satisfies these conditions then so does $g^{-1}$ with the same data. Suppose $G$ is contracting with nucleus $\mathcal{N}$ and that these conditions hold for some $g \in G$. The second condition implies that $g^l|_i = g^l$, so $g^l \in \mathcal{N}$ for all positive integers $l$. Since $\mathcal{N}$ is finite, we can take $l$ to be the smallest positive integer for which $g^l = 1$. Then $l = nl'$ for some integer $l'$, since otherwise $g^l$ has a nontrivial root permutation. But by the third condition, $g^{nl'}|_j = g^{ml'} = 1$, contrary to the minimality of $l$. 
\end{proof}

\begin{proof}[Proof of Proposition \ref{HTnoncontracting}]
Without loss of generality, we may assume that the two elements are $a = a_{01}$ and $b = a_{12}$. That is,
\[
a = (0\ 1)(1, 1, a, a, \ldots, a) \quad \quad b = (1\ 2)(b, 1, 1, b, \ldots, b).
\]
Then $ab = (0\ 1\ 2)(b, a, 1, ab, \ldots, ab)$ and $(ab)^3|_0 = ab$, $(ab)^3|_3 = (ab)^3$, so Lemma \ref{noncontractioncondition} applies with $g = ab$, $i = 3$, and $j = 0$.
\end{proof}

This means that any contracting self-similar group $(G, \mathsf{X}_k)$, $k > 3$, with generators $a_{ij}$ represents an uninteresting game, such as a 4-peg game with two legal moves: between pegs 0 and 1 and between pegs 2 and 3. In order to obtain modifications of the Hanoi Towers groups that are contracting, we are led to consider a larger class of possible generators.

\section{Hanoi Automorphisms and Groups} \label{HanoiAutomorphismsAndGroups}
In this section a class of automorphisms is defined that generalize the generators $a_{ij}$ that were used in the previous section.

\begin{definition}
An automorphism $a \in \Aut \mathsf{X}_k^{-\omega}$ is called a \emph{$k$-peg Hanoi automorphism} if there are disjoint subsets of $X_k,$ $P_a$ and $Q_a$, such that
\begin{enumerate}
\item $P_a \cup Q_a = \mathsf{X}_k$;
\item for each $i \in P_a$, $a|_i = 1$;
\item for each $j \in Q_a$, $a|_j = a$;
\item $\sigma_a$ fixes each element of $Q_a$.
\end{enumerate}
For $a \neq 1$, $P_a$ and $Q_a$, if they exist, are uniquely determined; we call them the sets of \emph{active} and \emph{inactive pegs} of $a$, respectively. The set of inactive pegs of the identity automorphism is defined to be $\mathsf{X}_k$.

Write $S_{k,q}$ for the set of automorphisms, $a$, of $X_k$ that have $q$ inactive pegs, and $S_k= \bigcup_q S_{k,q}$. The root permutation of an automorphism, $a$, is naturally regarded as a permutation of both $X_k$ and of $P_a$. 
\end{definition}

In terms of the game, $a \in S_{k, q}$ corresponds to ignoring some set $Q_a$ of $q$ pegs and moving the smallest disk among the remaining pegs, $P_a$, according to $\sigma_a$. For example, the 5-peg Hanoi automorphism $a = (0\ 1\ 2)(1, 1, 1, a, a) \in S_{5, 2} \subset S_5$ ignores pegs 3 and 4 and applies $(0\ 1\ 2)$ to the smallest disk among pegs 0, 1, and 2. A $k$-peg Hanoi automorphism $a$ is thus determined by $Q_a$ and $\sigma_a \in \Sym(\mathsf{X}_k \setminus Q_a)$. Note that $b = (0\ 1)(1, 1, 1, b, b) \in S_{5, 2} \subset S_5$ also has two inactive pegs 3 and 4, even though peg 2 is also fixed by $(0\ 1)$. If the smallest disk among pegs 0, 1, and 2 is on peg 2, then $b$ does nothing.

\begin{definition}
A \emph{$k$-peg Hanoi group} $(G, \mathsf{X}_k)$ is a group together with a fixed, state-closed generating set $S \subset S_k$ with $1 \in S$. The corresponding \emph{Hanoi game} is the $k$-peg game whose legal moves correspond to the non-identity generators. Moreover, the order of a Hanoi automorphism is the order of its root permutation; in particular this order is finite. Hence every Hanoi group is a self-similar group.
\end{definition}

For example, $H^{(k)}$ is the Hanoi group generated by $S_{k, k-2}$. The remainder of this section studies when Hanoi groups are contracting or non-contracting.

\begin{lemma}\label{lengthlemma}
Let $G$ be a Hanoi group. Consider $g \in G$ with representation $g = s_n \cdots s_2s_1$. Since the inverse of a Hanoi automorphism is another Hanoi automorphism (with the same inactive pegs and with the inverse root permutation), we can define $Q_i$ to be the set of inactive pegs of $s_i$. Define the \emph{essential set} of this representation to be $Q = \bigcap_{i = 1}^n Q_i$. If $g \neq 1$, then
\begin{enumerate}
\item $g|_j = g$ when $j \in Q$ and
\item $g = s_n \cdots s_2s_1$ is a minimal representation when for $j \notin Q$ $l(g|_j) < n$.
\end{enumerate}
\end{lemma}
\begin{proof}
Write $\sigma_i$ for the root permutation of $s_i$.
\begin{enumerate}
\item By Lemma \ref{sectionformula} for a given $j \in Q$ we have $g|_j = s_n|_j s_2|_j \cdots s_n|_j = g$.
\item Induct on $n$. Clearly $n = 1$ holds since $j \notin Q_1$ implies $s_1|_j = 1$. Assuming the claim for $n-1$, consider $g = s_n \cdots s_2s_1$, $n \ge 2$, and $j \notin Q$. If $j \notin Q_1$, then
\[
g|_j = (s_n \cdots s_2)|_{\sigma_1(j)}s_1|_j = (s_n \cdots s_2)|_{\sigma_1(j)}.
\]
By Lemma \ref{sectionformula}, this is a product of sections of the $s_i$, so $\l(g|_j) < n$. Now suppose $j \in Q_1$, so $g|_j = (s_n \cdots s_2)|_j s_1$. If $s_n \cdots s_2 = 1$, then $l(g|_j) = 1$. Otherwise, since $j \notin Q_2 \cap \cdots \cap Q_n$, we have $l((s_n \ldots s_2)|_j) \le n-2$ by the $n-1$ case, and so again $l(g|_j) < n$. \qedhere
\end{enumerate}
\end{proof}

\begin{proposition} \label{prenucleusprop}
The following conditions are equivalent for an element $g$ of a Hanoi group $(G, \mathsf{X}_k)$:
\begin{enumerate}
\item $g|_j = g$ for some $j \in \mathsf{X}_k$,
\item The essential set $Q$ is nonempty for some representation of $g$,
\item The essential set $Q$ is nonempty for all minimal representations of $g$.
\end{enumerate}
The set of all $g$ satisfying these conditions is called the \emph{prenucleus} of $G$, and is state-closed.
\end{proposition}

\begin{proof}
The three conditions clearly holds for the identity automorphism, so assume $g \neq 1$. Then 3 $\Rightarrow$ 2 is trivial, 2 $\Rightarrow$ 1 is part 1 of Lemma \ref{lengthlemma}, and the contrapositive of 1 $\Rightarrow$ 3 follows from part 2 of the same lemma.

It remains to show that the prenucleus is state-closed. Suppose $g = s_m \cdots s_2s_1$ satisfies the three conditions, and $j \in X_k$. Taking the section of $g$ at $j$ gives $$g|_j = s_m|_{j_m} \cdots s_2|_{j_2}s_1|_{j_1}$$ with $j_1=j$ (see Lemma \ref{sectionformula} for the notation). The sections of the generators are equal to either the identity automorphism or themselves and so the essential set of the representation of $g|_j$ contains the essential set of the representation of $g$ and is thus non-empty. Hence $g|_j$ is in the prenucleus for any $j \in X_k$ and $g$ in the prenucleus.
\end{proof}

\begin{theorem}\label{hanoicontraction}
A Hanoi group $G$ is contracting if and only if its prenucleus, $\mathcal{N}$, is finite. In this case, $\mathcal{N}$ is the nucleus of $G$.
\end{theorem}

\begin{proof}
If $G$ is contracting, its nucleus contains $\mathcal{N}$ by characterization 1 in Proposition \ref{prenucleusprop}. So if $\mathcal{N}$ is infinite, $G$ cannot be contracting.

Suppose $\mathcal{N}$ is finite, and consider $g \in G$ of length $n$. We show that any length-$n$ section $g|_{x_nx_{n-1}\ldots x_1}$ is in $\mathcal{N}$. Assume otherwise, then since $\mathcal{N}$ is state-closed, none of the sections $g|_{x_n}$, $g|_{x_nx_{n-1}}$, $\ldots$, $g|_{x_nx_{n-1}\ldots x_1}$ can be in $\mathcal{N}$. By characterization 3 in Proposition \ref{prenucleusprop}, each of these has a minimal representation with empty essential set. Thus by part 2 of Lemma \ref{lengthlemma}, the length decreases from one section to the next, and so $l(g|_{x_nx_{n-1}\ldots x_1}) \le 0$, a contradiction since $1 \in \mathcal{N}$ and $\mathcal{N}$ must be finite.

Thus $\mathcal{N}$ meets the condition for the finite set in the definition of the contraction property. We already noted that the nucleus contains $\mathcal{N}$. Since the nucleus is the smallest such set, it must equal $\mathcal{N}$.
\end{proof}

\begin{definition}
Let $G$ be a Hanoi group with prenucleus $\mathcal{N}$. By the second characterization in Proposition \ref{prenucleusprop}, every $g \in \mathcal{N}$ has a representation with a nonempty essential set. Define $d(g)$ as the least number of distinct generators among all such representations of $g$.
\end{definition}

\begin{definition}
Let $(G, \mathsf{X}_k)$ be a $k$-peg Hanoi group generated by $S$. For a subset $T \subset S$ with fixed indexing $T = \{ s_i \}_{i = 1}^M$, write $Q_i$ for the set of inactive pegs of $s_i$, and $\sigma_i$ for the root permutation of $s_i$. The \emph{essential set} of $T$ is defined to be $Q = \bigcap_{i = 1}^M Q_i$.

The subgroup of $\Sym(\mathsf{X}_k)$ generated by $\{ \sigma_i \}_{i = 1}^M$ acts by permutation on $\mathsf{X}_k$. Write $\Fix(T)$ for the set of elements of $\mathsf{X}_k$ that are fixed by every $\sigma_i$, and $\Orb_T(j)$ for the orbit of $j \in \mathsf{X}_k$ under this action.
\end{definition}

Note that the essential set of a representation (see Lemma \ref{lengthlemma}) coincides with the essential set of the set of generators it uses.

\begin{theorem} \label{hanoicontraction2}
Let $(G, \mathsf{X}_k)$ be a $k$-peg Hanoi group generated by $S$. Then $G$ is contracting if and only if the following condition holds:
\begin{description}
\item[($\ast$)] For any $T = \{ s_i \}_{i = 1}^M \subset S$ with nonempty essential set and any $j \notin \Fix(T)$, there exists some $i$ such that $\Orb_T(j) \cap Q_i = \emptyset$.
\end{description}
\end{theorem}

We use this lemma in the proof of the theorem.

\begin{lemma} \label{distinctgenlemma}
Assume a Hanoi group $G$ satisfies condition ($\ast$). Suppose $g \in G$ has a representation $g = a_m \ldots a_1$, with a nonempty essential set $Q$, using $M$ distinct generators $T = \{ s_i \}_{i = 1}^M$. Then for any $j \in \mathsf{X}_k$, either $g|_j = g$ or $d(g|_j) < M$.
\end{lemma}
\begin{proof}
If $j \in \Fix(T)$, then $g|_j = a_m|_j \cdots a_2|_j a_1|_j$ by Lemma \ref{sectionformula}. If $s_i|_j = 1$ for any $i$ then this representation does not use $s_i$, so $d(g|_j) < N$. Otherwise, $s_i|_j = s_i$ for every $i$ and $g|_j = g$.

Now suppose $j \notin \Fix(T)$. By Lemma \ref{sectionformula}, $g|_j = a_m|_{j_m} \cdots a_2|_{j_2} a_1|_{j_1}$, where $j_l \in \Orb_T(j)$ for all $l$, $m \ge l \ge 1$. Using $(\ast)$ to choose $i$ so that $\Orb_T(j)$ Again the essential set is nonempty, so $d(g|_j) < M$.
\end{proof}

\begin{proof}[Proof of Theorem \ref{hanoicontraction2}]
First suppose condition ($\ast$) holds. Let $\mathcal{N}$ be the prenucleus of $G$, and $\mathcal{N}_M = \{g \in \mathcal{N}\ |\ d(g) \le M\}$. We induct on $M$ to show that every $\mathcal{N}_M$ is finite. Clearly, $\mathcal{N}_1$ is finite. For $g \in \mathcal{N}_M$, $M \ge 2$, by Lemma \ref{distinctgenlemma}, every section of $g$ either is $g$ or is in $\mathcal{N}_{M-1}$. Because wreath recursion defines $g$ recursively, $g$ depends only on $\sigma_g$ and its sections that are not $g$. If $\mathcal{N}_{M-1}$ is finite, then there are only finitely many such distinguishing choices, and so $\mathcal{N}_M$ is also finite. Hence $\mathcal{N} = \mathcal{N}_{|S|}$ is finite, and $G$ is contracting by Theorem \ref{hanoicontraction}.

Now suppose ($\ast$) does not hold. That is, for some subset $T = \{s_i\}_{i = 1}^N$ of $S$ with nonempty essential set $Q$ and some $j \notin \Fix(T)$, we have $\Orb_T(j) \cap Q_i \neq \emptyset$ for every $i$, $1 \le i \le M$. It suffices to show that $H = \langle T \rangle \le G$ is infinite; indeed, since the prenucleus of $G$ contains $H$ by characterization 1 in Proposition \ref{prenucleusprop}, this would imply by Theorem \ref{hanoicontraction} that $G$ is not contracting.

We follow the idea in the proof of infinite cardinality of the Grigorchuk group (original proof in Russian in \cite{Gri80}; see Theorem 1.6.1 in \cite{Nek05} for an exposition in English). With $j \in \mathsf{X}_k$ chosen as above, let $H_j \le H$ be the subgroup of automorphisms whose root permutation fixes $j$. Then $\phi: H_j \rightarrow H$ defined by $\phi(g) = g|_j$ is a homomorphism. Since $j \notin \Fix(T)$, $H_j$ is a proper subgroup of $H$. For each generator $s_i$ of $H$, we construct below some $h_i \in H_j$ such that $\phi(h_i) = h_i|_j = s_i$, which shows that $\phi$ is a surjection. Then $H$ must be infinite since a proper subgroup $H_j < H$ maps onto $H$. 

Fix $i$, $1 \le i \le M$. By the assumption, $\Orb_T(j) \cap Q_i \neq \emptyset$, so there exists some $s_{n_L} \ldots s_{n_2}s_{n_1} \in H$ with root permutation satisfying $\sigma_{n_L} \ldots \sigma_{n_2}\sigma_{n_1}(j) \in Q_i$. Let $j_1 = j$ and, for $l$, $L \ge l \ge 1$, let $j_{l + 1} = \sigma_{n_{l}}(j_{l})$. 
Let $j_1=j$ and let $j_{l+1} = \sigma_{n_{l}}(j_{l})$ for $1 \le l \le L$.
This implies $s_{n_l}|_{j_l} = s_{n_l|_{j_{l+1}}} = 1$ for all $l$, $L \ge l \ge 1$.

Take
\begin{align*}
h_i &= (s_{n_L} \cdots s_{n_2}s_{n_1})^{-1}s_i(s_{n_L} \cdots s_{n_2}s_{n_1}) \\
&= (s_{n_1}^{-1}s_{n_2}^{-1} \cdots s_{n_L}^{-1}) s_i (s_{n_L} \cdots s_{n_2} s_{n_1})
\end{align*}
with root permutation $(\sigma_{n_1}^{-1}\sigma_{n_2}^{-1} \cdots \sigma_{n_L}^{-1}) \sigma_i (\sigma_{n_L} \cdots \sigma_{n_2}\sigma_{n_1})$. Since $\sigma_i$ fixes \\ $\sigma_{n_L} \cdots \sigma_{n_2}\sigma_{n_1}(j) = j_{L + 1} \in Q_i$,
\begin{align*}
\sigma_{n_l} \cdots \sigma_{n_2}\sigma_{n_1}(j) &= j_{l + 1} \mbox{ for } L \ge l \ge 1 \\
\sigma_i \sigma_{n_L} \cdots \sigma_{n_2}\sigma_{n_1}(j) &= j_{L + 1} \\
\sigma_{l}\sigma_{n_{l+1}} \cdots \sigma_{n_L} \sigma_i \sigma_{n_L} \cdots \sigma_{n_2}\sigma_{n_1}(j) &= j_l \mbox{ for } 2 \le l \le L \\
\sigma_{n_1}\sigma_{n_2} \cdots \sigma_{n_L} \sigma_i \sigma_{n_L} \cdots \sigma_{n_2}\sigma_{n_1}(j) &= j_{l - 1}.
\end{align*}
The last equation shows $h_i \in H_j$. Using the remaining equations in Lemma \ref{sectionformula},
\[
h_i|_j = \left( s_{n_1}^{-1}|_{j_2} s_{n_2}^{-1}|_{j_3} \cdots s_{n_l}^{-1}|_{j_{l+1}} \right) s_i|_{j_{l+1}} \left( s_{n_l}|_{j_l} \cdots s_{n_2}|_{j_2}s_{n_1}|_{j_1} \right) = s_i,
\]
as we wanted.
\end{proof}

Given $S \subset S_k$, checking ($\ast$) is a finite computation. Moreover, since $T_1 \subset T_2 \subset S$ implies $\Fix(T_2) \subset \Fix(T_1)$ and $\Orb_{T_1}(j) \subset \Orb_{T_2}(j)$, it suffices to check the maximal sets among subsets $T \subset S$ with nonempty essential set.

\section{Contracting Hanoi Groups  and Their Limit Spaces} \label{Hck}

In view of Theorems \ref{hanoicontraction} and \ref{hanoicontraction2}, $H^{(k)}$ is noncontracting for $k \ge 4$ because two generators had distinct but overlapping sets of inactive pegs (common inactive peg 3 in the proof of Proposition \ref{HTnoncontracting}). This is avoided if only one peg is inactive, which leads us to the following definition.

\begin{definition}
$H_c^{(k)}$ is the $k$-peg Hanoi group generated by $S_{k, 1}$.
\end{definition}

\begin{corollary} \label{modifiedhanoicontracting}
$H_c^{(k)}$ is contracting for every $k \ge 3$.
\end{corollary}
\begin{proof}
Since the prenucleus of $H_c^{(k)}$ is $S_{k, 1}$, the result is immediate from Theorem \ref{hanoicontraction}.
\end{proof}

Alternatively, this follows since $H_c^{(k)}$ is contained in the group of \emph{bounded automorphisms} of $\mathsf{X}_k^*$ (see \cite{BN03} and Theorem 3.9.12 in \cite{Nek05}). In Section \ref{Symmetry}, we will use Theorem \ref{hanoicontraction2} in a case where this result is not applicable.

For $k = 3$, $S_{3, 1} = \{ 1, a_{01}, a_{12}, a_{02} \}$, and $H^{(3)} = H_c^{(3)}$. We can thus view $H_c^{(k)}$ as higher-peg analogues that preserve contraction. In the Hanoi game corresponding to $H_c^{(k)}$, the $n$-th smallest disk can be moved only when the $n-1$ smallest disks are on the same peg. As a result, for any number of pegs, the optimal move count $h_n$ for moving $n$ disks from one peg to another satisfies the same recurrence $h_n = 2h_{n-1} + 1$ as for three pegs (see \cite{BH99}). Additional pegs, while providing more solution paths, do not make the game any shorter.

Because these groups are contracting we have a well-defined limit space to analyze as was done in \cite{GS08} for $H^{(3)}.$ By construction $S_{k,1}$ is the nucleus of $H_c^{(k)}$ so by Proposition \ref{nucleusasymtotic} asymptotically equivalent pairs of left-infinite sequences have the form $(\ldots llliv, \ldots llljv)$, where $i, j, l \in \mathsf{X}_k$ are all distinct and $v \in \mathsf{X}_k^*$. For simplicity's sake abbreviate $\mathcal{J}_{H_c^{(k)}}$ as $\mathcal{J}^{(k)}$ for the limit space of $H_c^{(k)}$. We will show for each $k \ge 3$ that $\mathcal{J}^{(k)}$ can be obtained as a self-similar set $K^{(k)}$ in $\mathbb{R}^{k-1}$. We first, however, need some basic definitions and results about self-similar sets.

\begin{definition}
(Definitions 1.2.1 and 1.3.1 in \cite{Kig01}) For a finite set $S = \{ s_i \}_{i=1}^{N}$, the \emph{shift space} $\Sigma(S) = S^{-\omega}$ is the set of all left-infinite sequences of elements of $S$, with the product topology where $S$ has the discrete topology. For each $s_i \in S$, define $s_i: \Sigma(S) \rightarrow \Sigma(S)$ by $s_i(\ldots x_3x_2x_1) = \ldots x_3x_2x_1s_i$.

Let $K$ be a compact metric space. For each $s_i \in S$, let $F_i$ be a continuous injection from $K$ to itself. Then $(K, S, \{ F_i \}_{s_i \in S})$ is called a \emph{self-similar structure} if there exists a continuous surjection $\pi: \Sigma(S) \rightarrow K$ such that $\pi \circ s_i = F_i \circ \pi$ for any $s_i \in S$. We call $K$ a \emph{self-similar set}.
\end{definition}

\begin{theorem}\label{selfsimilarstructure}
(Theorems 1.1.4 and 1.2.3 in \cite{Kig01}) Let $(X, d)$ be a complete metric space and $S = \{ s_1, s_2, \ldots, s_N \}$. For each $s_i \in S$, let $F_i: X \rightarrow X$ be a contraction. For any $A \subset X$, define $F(A) = \bigcup_{i = 1}^N F_i(A)$. Then there exists a unique nonempty compact set $K$ satisfying $F(K) = K$.

For $w = x_m \ldots x_2x_1$, $x_i \in S$, set $F_w = F_{x_m} \circ \cdots \circ F_{x_2} \circ F_{x_1}$ and $K_w = F_w(K)$. Each $K_w$, $|w| = m$, is called an \emph{$m$-cell}. Then for any $w = \ldots x_2x_1 \in \Sigma(S)$, $\bigcap_{m \ge 1}K_{x_m \ldots x_2x_1}$ contains only one point. The map $\pi: \Sigma(S) \rightarrow K$ defined by taking $\pi(w)$ as this unique point is a continuous surjection. Moreover, $\pi \circ s_i = F_i \circ \pi$ for all $s_i \in S$.
\end{theorem}

Fix $k \ge 3$. Let $\{ p_i\ |\ i \in \mathsf{X}_k \}$ be the vertices of a regular $k$-simplex in $\mathbb R^{k-1}$. For each $i$, the points $\{ p_j\ |\ j \neq i \}$ form the vertices of a regular $(k-1)$-simplex; let $q_i$ be its centroid. Define $F_i$, $i \in \mathsf{X}_k$, as the unique affine map satisfying
\[
F_i(p_j) = \left\{ \begin{array}{ll}
             p_j & \mbox{if $i = j$} \\
             q_j & \mbox{otherwise}.\end{array} \right.
\]
Since each $F_i$ is a contraction with respect to the Euclidean metric, by Theorem \ref{selfsimilarstructure} there is a unique nonempty compact set $K^{(k)}$ with $F(K^{(k)}) = K^{(k)}$, and $( K^{(k)}, \mathsf{X}_k, \{ F_i \}_{i \in \mathsf{X}_k} )$ is a self-similar structure.

Before proving that $K^{(k)}$ is homeomorphic to $\mathcal{J}^{(k)}$, we need an observation.

\begin{figure}[htbp] 
\begin{center}
\includegraphics[width=50mm]{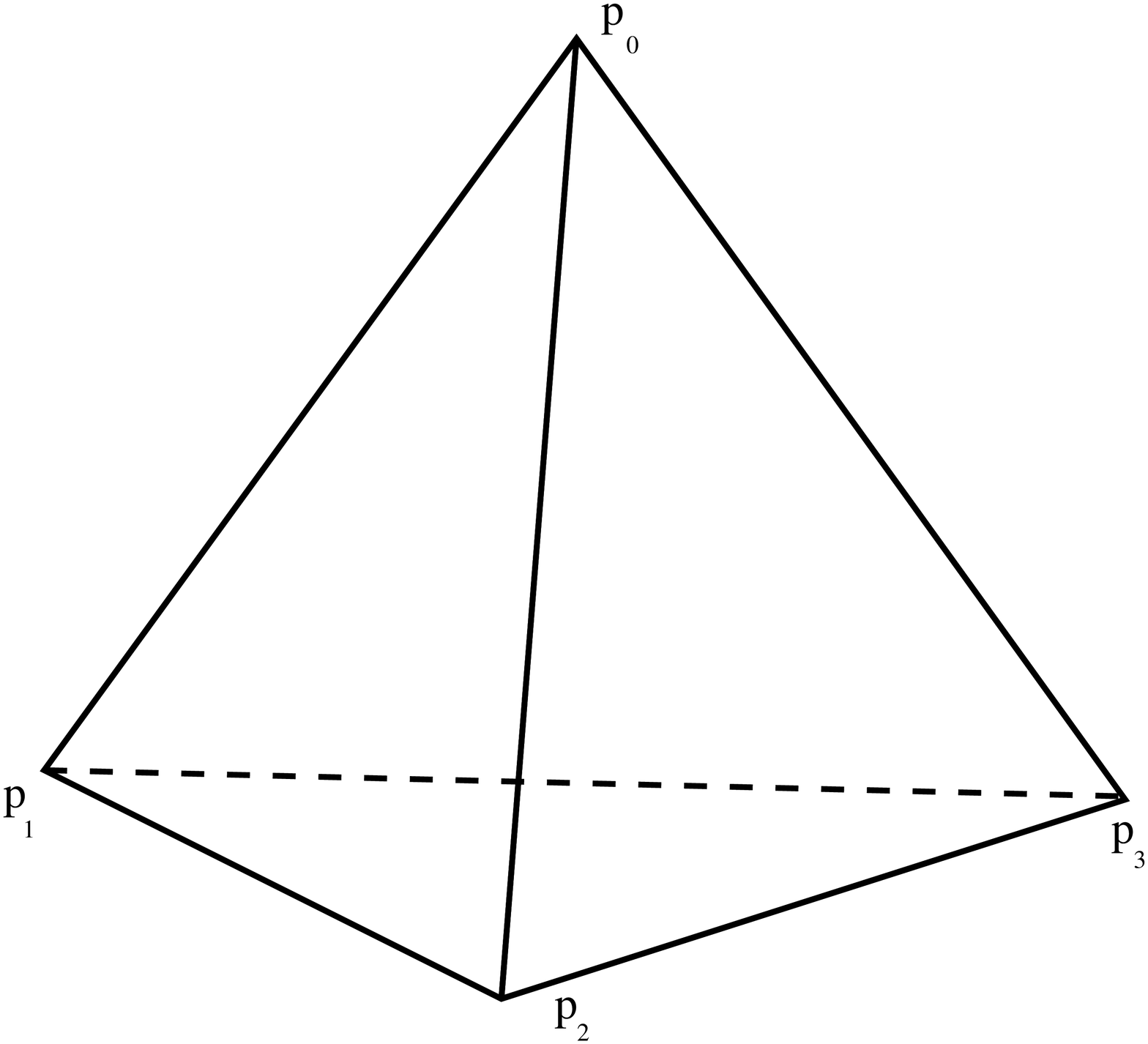} \includegraphics[width=50mm]{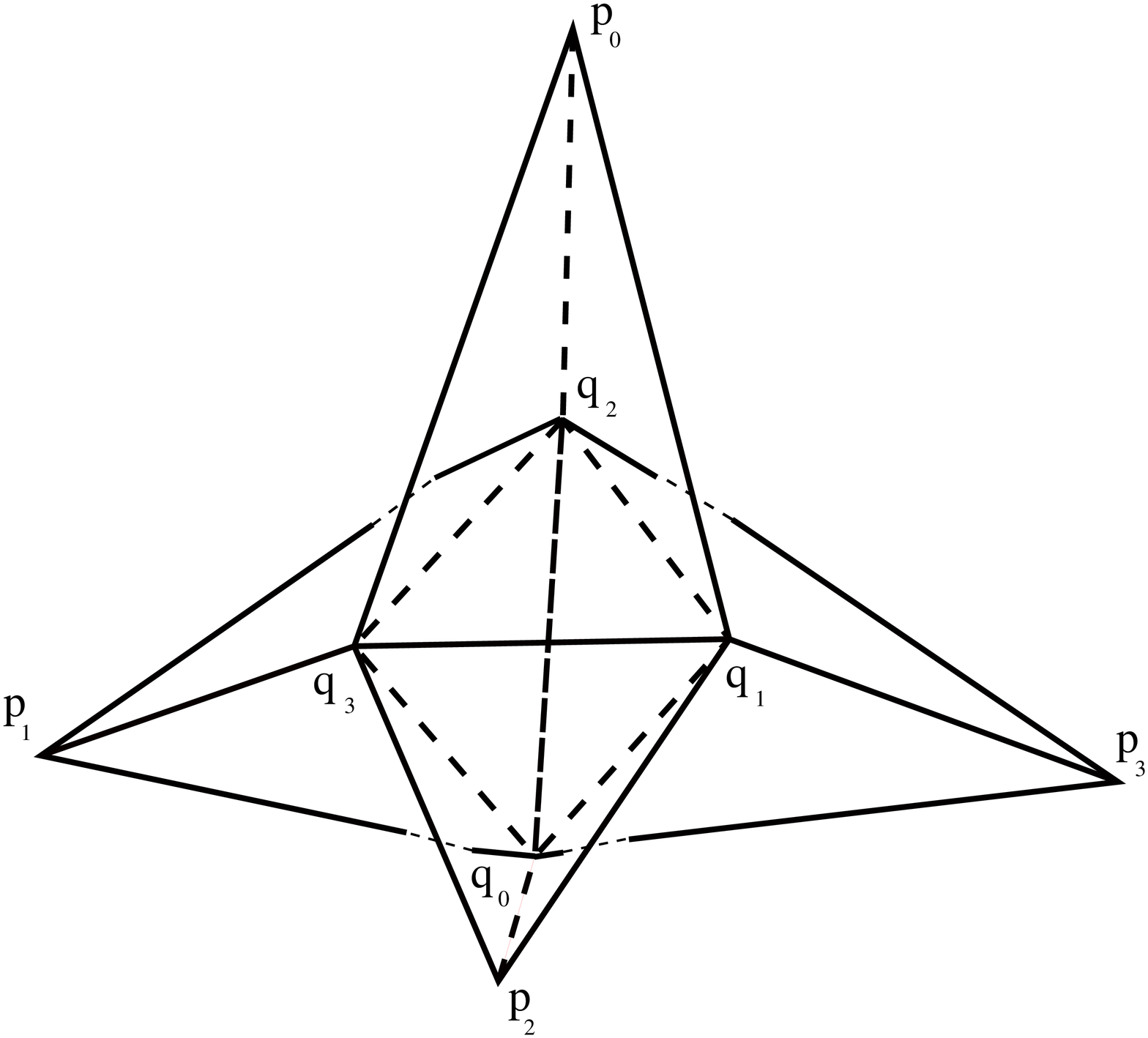}
\end{center}
\caption{$E^{(4)}$ and $F(E^{(4)})$ in $\mathbb{R}^3$}
\label{k=4example}
\end{figure}

\begin{lemma} \label{E'intersection}
Let $E^{(k)}$ be the $k$-simplex with vertices $\{p_i\ |\ 0 \le i \le k-1 \}$. Then for $i \neq j$,
\[
F_i(F(E^{(k)})) \cap F_j(F(E^{(k)})) = \{ q_l\ |\ l \neq i, j \}.
\]
Since $K \subset F(E^{(k)})$, it follows that $F_i(K) \cap F_j(K) = \{ q_l\ |\ l \neq i, j \}$.
\end{lemma}

\begin{proof}
Without loss of generality, assume $i = 0$ and $j = 1$. It is clear that $F_0(F(E^{(k)})) \cap F_1(F(E^{(k)})) \supset \{ q_l\ |\ l \ge 2 \}$, so we show the reverse inclusion.

The \emph{convex hull} of a finite set of points $X = \{ x_i \}_{i = 1}^N$, denoted $\Conv(X)$, is the set of points $\sum_{i = 1}^N \alpha_ix_i$, where $\alpha_i$ are nonnegative reals satisfying $\sum_{i = 1}^N \alpha_i = 1$. Given such a representation, we call $\alpha_i$ the \emph{coefficient} of $x_i$

 $E = \Conv(\{p_i\ |\ 0 \le i \le k-1 \} )$. For a simplex, it is easily shown that $(\alpha_0, \ldots, \alpha_{k-1})$ gives a unique representation of every point in $E$. Since the $F_i$ are affine maps,
\begin{align*}
f_0(E) &= \Conv(\{p_0, q_1\} \cup \{q_i\ |\ i \ge 2\}) \\
f_1(E) &= \Conv(\{p_1, q_0\} \cup \{q_i\ |\ i \ge 2\}).
\end{align*}
Since $q_j$ is a weighted average of $\{p_i\ |\ i \neq j\}$, $\alpha_i > 0$ for every $i \neq j$. Thus for points of $F_0(E)$, $\alpha_0 > \alpha_1$ if the coefficient of either $p_0$ or $q_1$ is nonzero. Similarly, for points of $F_1(E)$, $\alpha_1 > \alpha_0$ if the coefficient of either $p_1$ or $q_0$ is nonzero. Hence $F_0(E) \cap F_1(E) \subset \Conv(\{q_i\ |\ i \ge 2\})$, and since $F(E^{(k)}) \subset E$,
\begin{eqnarray} \label{lemmaeq}
F_0(F(E^{(k)})) \cap F_1(F(E^{(k)})) \subset \Conv(\{q_i\ |\ i \ge 2\}).
\end{eqnarray}
For any $q_j$, at least one of $\alpha_0$ and $\alpha_1$ is positive. Then $p_i$ is the only element of $F_i(F(E))$ that can possibly lie in $\Conv(\{p_j\ |\ j \ge 2\})$. Since $F(E^{(k)}) = \bigcup F_i(E)$, we have $F(E^{(k)}) \cap \Conv(\{p_j\ |\ j \ge 2\}) \subset \{ p_i\ |\ j \ge 2 \}$. Applying $F_i$, $i \in \{0, 1\}$, to both sides and taking their intersection,
\[
[F_0(F(E^{(k)})) \cap F_1(F(E^{(k)}))] \cap \Conv(\{ q_j\ |\ j \ge 2 \}) = \{ q_j\ |\ j \ge 2 \}.
\]
The result follows by this and (\ref{lemmaeq}).
\end{proof}

\begin{proposition}
For all $k \ge 3$, $K^{(k)}$ is homeomorphic to $\mathcal{J}^{(k)}$.
\end{proposition}

\begin{proof}
Let $\pi: \mathsf{X}_k^{-\omega} \rightarrow K^{(k)}$ be as in Proposition \ref{selfsimilarstructure}, and suppose $\pi(\ldots x_3x_2i) = \pi(\ldots y_3y_2j) = p$ with $i \neq j$. Then $p \in F_i(K) \cap F_j(K) \subset \{ q_l\ |\ l \neq i, j \}$ by Lemma \ref{E'intersection}. It is easy to see that $p = q_l$ if and only if $x_n = y_n = l$ for all $n \ge 2$. Since the $F_i$ are injective, by Proposition 1.2.5 in \cite{Kig01}, $\pi(w) = \pi(\tau)$ for $w \neq \tau$ if and only if $(w, \tau) = (\ldots lliv, \ldots lljv)$ for $i, j, l \in \mathsf{X}_k$ all distinct and $v \in \mathsf{X}_k^*$.

Since this is the asymptotic equivalence relation in $\mathcal{J}^{(k)}$, $\pi$ induces a continuous bijection $p: \mathcal{J}^{(k)} \rightarrow K^{(k)}$. The limit space,  $\mathcal{J}^{(k)}$, is compact by Theorem 3.6.3 in \cite{Nek05}, and $K^{(k)}$ is Hausdorff; hence $p$ is a homeomorphism.
\end{proof}

It is known that $K^{(3)}$ is the Sierpi\'{n}ski gasket. Unlike in the usual construction, each $F_i$ here involves a reflection as well as a contraction. This reflects some of the self-similarity of the Hanoi Towers game (see Section \ref{HanoiTowerAutomaton}).

Except for when $k = 3$, $F_i$ is not a similitude since $F_i(E^{(k)})$ is not a regular $k$-simplex; that is, there is no constant $r$ for which $|F_i(x) - F_i(y)| = r|x - y|$. While the $F_i$ become similitudes if $q_i$ is taken outside $E^{(k)}$, the resulting invariant set is not homeomorphic to $\mathcal{J}^{(k)}$. Since the available theory is limited for sets invariant under affine contractions, the Euclidean metric on $K^{(k)}$ is of limited use. However, analytic structures such as energy, measure, and Laplacians have been constructed on the following special class of self-similar structures, independent of an Euclidean embedding \cite{Kig01}.

\begin{definition}
For a self-similar structure $(K, S, \{F_s\}_{s \in S})$, define the \emph{critical set} $\mathcal{C}$ and the \emph{post critical set} $\mathcal{P}$ by $\mathcal{C} = \pi^{-1}(\cup_{i, j \in S, i \neq j} (F_i(K) \cap F_j(K)))$ and $\mathcal{P} = \cup_{n \ge 1}\sigma^n(\mathcal{C}_{\mathcal{L}})$. If $\mathcal{P}$ is finite, the self-similar structure is said to be \emph{post critically finite} (\emph{p.c.f.} for short).
\end{definition}

\begin{proposition}
The self-similar structure $(K^{(k)}, \mathsf{X}_k, \{ F_i \}_{i \in \mathsf{X}_k})$ is p.c.f.
\end{proposition}
\begin{proof}
By Lemma \ref{E'intersection}, $\mathcal{C} = \{ \ldots lli\ |\ i \neq l\}$ and hence $\mathcal{P} = \{ \ldots ll\ |\ l \in \mathsf{X}_k \}$.
\end{proof}

Analysis can thus be developed on $K^{(k)}$ following the standard theory of p.c.f. self-similar sets (see \cite{Kig01}), and moreover, we expect the high symmetry of $K^{(k)}$ to simplify computations. As examples, in the rest of this section we compute the Hausdorff and spectral dimensions of $K^{(k)}$. We provide only a sketch of the theory; for details see \cite{Kig01} and Chapter 4 of \cite{Str06}.

Let $V_0 = \pi(\mathcal{P}) = \{ p_i\ |\ i \in \mathsf{X}_k \}$, and let $\Gamma_0$ be the complete graph on $V_0$. Inductively define $V_{m+1} = F(V_m)$ and $\Gamma_{m+1}$, a graph on $V_{m+1}$, where $x$ and $y$ are connected in $\Gamma_{m+1}$ if $x = F_i(x')$, $y = F_i(y')$ for some $F_i$ and $x', y' \in V_m$ connected in $\Gamma_m$. Given positive weights $\{ c_{ij}\ |\ i, j \in \mathsf{X}_k, i < j \}$, define an energy form on $V_0$ by
\[
\mathcal{E}_0(u) = \sum_{i < j}c_{ij}(u(p_i) - u(p_j))^2
\]
for any function $u$ on $V_0$. We wish to construct graph energies $\mathcal{E}_m$ on $\Gamma_m$ subject to two conditions.
\begin{description}
\item[(E1)] For some fixed \emph{renormalization factors} $\{ r_i \}_{i = 0}^{k-1}$,
\[
\mathcal{E}_m(u) = \sum_{i = 0}^{k-1}r_i^{-1}\mathcal{E}_{m-1}(u \circ F_i).
\]
\item[(E2)] Given $u$ on $V_{m-1}$, then an extension $u'$ of $u$ to $V_m$ is a function on $V_m$ such that $u'|_{V_{m-1}} = u$. Let $\tilde{u}$ be the extension of $u$ that minimizes $\mathcal{E}_m$. Then
\[
\mathcal{E}_m(\tilde{u}) = \mathcal{E}_{m-1}(u).
\]
\end{description}
Finding $\{ c_{ij} \}$ and $\{ r_i \}$ satisfying these conditions for a given self-similar structure, called the \emph{renormalization problem}, is in general highly nontrivial. Below, we find a solution for $K^{(k)}$ that is invariant under the full symmetry of $V_0$.

\begin{theorem} \label{renormalization}
The weights $c_{ij} = 1$, $i < j$, and renormalization factor $r_i = r = \frac{k(k - 2)}{k^2 - k - 1}$ are a solution to the renormalization problem for $(K^{(k)}, \mathsf{X}_k, \{ F_i \}_{i \in \mathsf{X}_k})$.
\end{theorem}

\begin{proof}
Let $p_i \in V_0$, $q_i \in V_1 \setminus V_0$. Take $c_{ij} = 1$ and, for now, $r_i = 1$. For a function $u$ on $V_0$ define $u(p_i) = a_i$,
\[
\mathcal{E}_0(u) = \sum_{i < j}(a_i - a_j)^2.
\]
For a function $u$ on $V_0$ defined by $u(p_i)=a_i$, an extension $u'$ of $u$ to $V_1$ is determined by its values on $V_1 \setminus V_0 = \{ q_i \}$; let $u'(q_i) = x_i$. Each edge of $\Gamma_1$ connecting $q_i$ and $p_j$, $i \neq j$, appears in one of the $k$ subgraphs $F_l(\Gamma_0)$ of $\Gamma_1$, while each edge connecting $q_i$ and $q_j$, $i \neq j$, appears in $k-2$ of them. Thus by (E1),
\[
\mathcal{E}_1(u') = (k-2)\sum_{i < j}(x_i - x_j)^2 + \sum_{i < j} (x_i - a_j)^2.
\]
The values of $x_i$ minimizing this determine $\tilde{u}$ on $V_m$. If $\mathcal{E}_1(\tilde{u})/\mathcal{E}_0(u) = \lambda$ for all $u$, then replacing $r_i$ by $r_i/\lambda$ makes (E2) hold for $m = 1$.

Setting  $\frac{1}{2}\frac{\partial\mathcal{E}_1}{\partial x_i} = 0$ yields
\begin{align*}
(k - 2)\sum_{j: j \neq i}(x_i - x_j) + \sum_{j: j \neq i}(x_i - a_j) &= 0 \\
[(k - 2)(k - 1) + (k - 1)]x_i &= (k - 2)\sum_{j: j \neq i}x_j + \sum_{j: j \neq i} a_j
\end{align*}
for each $i$. Summing $k$ such equations, one for each $i \in \mathsf{X}_k$, we get $(k-1)^2X = (k-1)(k-2)X + (k-1)X$ or $X = A$, where $X = \sum x_i$ and $A = \sum a_i$. So adding $(k-2)x_i$ to both sides to the previous equation gives
\begin{align*}
[(k-2)(k-1) + (k-1) + (k-2)]x_i &= (k-2)X + A - a_i \\
(k^2 - k - 1)x_i &= (k-1)X - a_i.
\end{align*}
Then $a_i - a_j = (k^2 - k - 1)(x_j - x_i)$, so
\[
\mathcal{E}_0(u) = \sum_{i < j}(a_i - a_j)^2 = (k^2 - k - 1)^2\sum_{i < j}(x_i - x_j)^2.
\]
Using $a_i = (k-1)X - (k^2 - k - 1)x_i$ to write $\mathcal{E}_1(\tilde{u})$ in terms of $k$ and symmetric sums $\sum_i x_i^2$ and $\sum_{i < j}x_ix_j$, we arrive at
\[
\mathcal{E}_1(\tilde{u}) = k(k - 2)(k^2 - k - 1) \sum_{i < j}(x_i - x_j)^2
\]
\[
\Rightarrow \lambda = \frac{\mathcal{E}_1(\tilde{u})}{\mathcal{E}_0(u)} = \frac{k^2 - k - 1}{k(k - 2)}.
\]
Hence (E2) holds for $m = 1$ with $r_i = \frac{k(k - 2)}{k^2 - k - 1}$.

For the extension from $\mathcal{E}_{m-1}$ to $\mathcal{E}_m$, note that each point in $V_m \setminus V_{m-1}$ lies in a unique $m$-cell, and that within each $m$-cell the minimization problem is what we just solved. Since $\mathcal{E}_m$ is the sum of the energy contribution from each cell, (E1) and (E2) hold for any $m$.
\end{proof}

Because $\mathcal{E}_m$ is nondecreasing by (E2), we define $\mathcal{E}(u) = \lim_{n \rightarrow \infty}\mathcal{E}_m(u)$, allowing values of $+\infty$. The \emph{effective resistance} $R$ is then defined on $V_* = \bigcup_{m \ge 0} V_m$ by
\[
R(x, y)^{-1} = \min\{ \mathcal{E}(u): u(x) = 0 \mbox{ and } u(y) = 1\}.
\]
This is independent of $m$. By construction $V_*$ is dense in $K$ and $R(x, y)$ is uniformly continuous in $x$ and $y$, so $R$ extends to $K \times K$ as a metric on $K$.

\begin{proposition}\label{diamcontractiveratio}
Let $R$ be the effective resistance metric on $K^{(k)}$, and $r = \frac{k(k - 2)}{k^2 - k - 1}$ as in Theorem \ref{renormalization}. There exist constants $c_1'$, $c_2'$ such that \\
(a) if $x$ and $y$ are in the same or adjacent $m$-cell, then $R(x, y) \le c_1'r^m$; \\
(b) otherwise $R(x, y) \ge c_2'r^m$.
\end{proposition}

\begin{proof}
Lemma 1.6.1 in \cite{Str06} proves this for the Sierpi\'{n}ski gasket. Section 4.4 of the same reference argues why this holds for any p.c.f. self-similar set.
\end{proof}

For self-similar Euclidean sets generated by contractive \emph{similitudes} of ratio $r_i$, the Hausdorff dimension $d_H$ under the Euclidean metric satisfies Moran's formula \cite{Mor46}, $\sum_{i = 0}^{k-1} r_i^{d_H} = 1$. The above Proposition says that the contraction maps have ratio roughly $r$. Here, we need a generalization from Kigami (Theorem 1.5.7 in \cite{Kig01}). We use a special case.

\begin{theorem} (Corollary 1.3 in \cite{Kig95}) \label{morangeneralization}
Let $(X, d)$ be a complete metric space, and let $K$ be the self-similar invariant set defined by contractions $\{f_i\}_{i \in S}$. Suppose there exist constants $0 < r < 1$, $0 < c_1, c_2$ and $M > 0$ such that \\
(1) for all $w \in \mathsf{X}^n$,
\[
\diam(K_w) \le c_1r^n;
\]
(2) for all $x \in K$ and $n \ge 0$,
\[
\#\{ w: w \in \mathsf{X}^n, d(x, K_w) \le c_2r^n\} \le M.
\]
Then the Hausdorff dimension of $K$ with respect to $d$ is $-\log{|S|}/\log{r}$. 
\end{theorem}

\begin{corollary}
The space $K^{(k)}$ has Hausdorff dimension $\frac{\log{k}}{\log{(k^2 - k - 1)}- \log{(k(k-2))}}$ with respect to the effective resistance metric.
\end{corollary}

\begin{proof}
Let $c_1'$ and $c_2'$ be as in Proposition \ref{diamcontractiveratio}, and take $c_1 = c_1'$, $c_2 = c_2'r$, $M = k + 1$. By (a) in the same proposition, $\diam(K_w) \le c_1r^n$ for $w \in \mathsf{X}^n$. Take any $x \in K = K^{(k)}$ and $n \ge 0$, and consider $w \in \mathsf{X}_k^n$. If $K_w$ is not one of the $M = k+1$ $n$-cells that either contain $x$ or is adjacent to the $n$-cell containing $x$, then by (b), $R(x, K_w) \ge c_2'r^n > c_2r^n$. Hence Theorem \ref{morangeneralization} applies.
\end{proof}

The renormalization problem is often stated through the terminology of \emph{harmonic structure}. Write $l(V_0)$ for the space of functions on $V_0$. The quadratic form $-\mathcal{E}_0$ can be written $-\mathcal{E}_0(u) = u^TDu$ for some symmetric linear operator $D: l(V_0) \rightarrow l(V_0)$. For $\mathcal{E}_0$, constructed above, $D$ has the $k$-by-$k$ matrix representation
\[
D =
\begin{pmatrix}
-(k-1) & 1 & \cdots & 1 \\
1 & -(k-1) & \cdots & 1 \\
\vdots  & \vdots  & \ddots & \vdots  \\
1 & 1 & \cdots & -(k-1)
\end{pmatrix}
\]
with respect to the standard function basis.  Given $\{r_i\}_{i \in \mathsf{X}_k}$, define linear operators $H_m: l(V_m) \rightarrow l(V_m)$ inductively by $H_0 = D$ and 
\[
H_m = \sum_{i \in \mathsf{X}_k} r_i^{-1}R_i^{-1}H_{m-1}R_i,
\]
where $R_i: l(V_m) \rightarrow l(V_{m-1})$ is the natural restriction defined by $R_i(u) = u \circ F_i$. Then condition (E2) for energy together with $-\mathcal{E}_0 = u^TDu$ implies that $-\mathcal{E}_m(u) = u^T H_m u$ for all $m \ge 0$. Decompose $H_1$ as 
\[
H_1 = \left(
\begin{array}{c|c}
T & J^T \\ \hline
J & X
\end{array}
\right),
\]
where $T$ acts on $l(V_0)$ and $X$ acts on $l(V_1 \setminus V_0)$. In the terminology of \cite{KL93}, $(D, \{r_i\}_{i \in \mathsf{X}_k})$, with suitable additional conditions on $D$, is called a \emph{harmonic structure} on $(K^{(k)}, \mathsf{X}_k, \{ F_i \}_{i \in \mathsf{X}_k})$ if
\begin{eqnarray} \label{renormeq}
D = \lambda(T - J^TX^{-1}J)
\end{eqnarray}
for some constant $\lambda$. A harmonic structure is said to be \emph{regular} if $\lambda > r_i$ for all $i \in \mathsf{X}_k$.

\begin{corollary}
The pair $(D, \{r, \ldots, r\})$ for $D$ defined above and $r = \frac{k(k-2)}{k^2 - k - 1}$ is a \emph{regular harmonic structure} on $(K^{(k)}, \mathsf{X}_k, \{ F_i \}_{i \in \mathsf{X}_k})$ with $\lambda = 1$.
\end{corollary}

\begin{proof}
Minimizing $\mathcal{E}_1(\tilde{u}) = \tilde{u}^T H_1 \tilde{u}$ as in the proof of Theorem \ref{renormalization} shows that (\ref{renormeq}) is equivalent to $\mathcal{E}_1(\tilde{u}) = \lambda\mathcal{E}_0(u)$. Recall that in Theorem \ref{renormalization} we replaced $r_i$ by $r_i/\lambda$ to satisfy condition (E2), which is just the above condition with $\lambda = 1$. Here we instead leave $\lambda$ apart from $\{r_i\}$. Therefore, finding $D$, $\{r_i\}$, and $\lambda$ under these conditions is an equivalent reformulation of the renormalization problem, differing merely in the choice of the fundamental analytic structure (harmonic structure instead of a self-similar energy form). Hence the result of Theorem \ref{renormalization} carries over. The harmonic structure is regular since $\lambda = 1 > r_i$.
\end{proof}

By further specifying the standard Bernoulli (self-similar) measure on $K^{(k)}$, the standard Laplacian on $K^{(k)}$ can be defined as the limit of the operators $H_m$ on $\Gamma_m$. The \emph{spectral dimension} of the harmonic structure describes the spectral asymptotics of this Laplacian.

\begin{theorem}
The \emph{spectral dimension} $d_S$ of $K^{(k)}$ with the harmonic structure $(D, \{r, \ldots, r\})$ defined above is $\frac{2\log{k}}{ \log(k^2 - k - 1) - \log(k - 2)}$.
\end{theorem}
\begin{proof}
By Theorem A.2 in \cite{KL93},
\[
\sum_{i = 0}^{k-1}\gamma_i^{d_S} = 1,
\]
where $\gamma_i = \sqrt{\frac{r_i\mu_i}{\lambda}}$, with $\mu_i = \frac{1}{k}$ for the standard Bernoulli measure. We omit the computation.
\end{proof}

As a remark, since this self-similar structure is fully symmetric on $V_0$, by Proposition 3.1 in \cite{Betal08}, the \emph{spectral decimation} method allows the explicit computation of the spectrum of the associated Laplacian to this energy form with multiplicities.

\section{Symmetric Contracting Hanoi Groups} \label{Symmetry}
The pegs of the original Hanoi Towers game are indistinguishable in the sense that the same type of legal move exists between any two pegs. By placing symmetry conditions on the Hanoi groups, we obtain corresponding Hanoi games that preserve this aspect of the game.

The permutation group on $\mathsf{X}_k$, $\Sym(\mathsf{X}_k)$, has a natural action on $S_k$; for $\phi\in\Sym(\mathsf{X}_k)$ and $a \in S_k$, define $\phi \cdot a$ as the Hanoi automorphism with inactive pegs $\phi(Q_a)$ and root permutation $\phi\sigma_a\phi^{-1}$ on the active pegs $\phi(\mathsf{X}_k \setminus Q_a)$. In the game, this corresponds to rerepresenting the same move by relabeling each peg $i$ as $\phi(i)$.

\begin{definition}
A $k$-peg Hanoi group $G$ with generating set $S$ is said to be \emph{fully symmetric} if $S$ is closed under $\Sym(\mathsf{X}_k)$; that is, $\phi \cdot S \subset S$ for all $\phi \in \Sym(\mathsf{X}_k)$.
\end{definition}

In fact, since $S$ is finite and the action faithful, $\phi \cdot S = S$. More concretely, if $g \in S$, then $S$ contains every $k$-peg Hanoi automorphism with the same number of inactive pegs and the same cycle type of the root permutation as $g$.

\begin{theorem}\label{uniquecontracting}
For any $k \ge 3$, every fully symmetric and contracting $k$-peg Hanoi group is a subgroup of $H_c^{(k)}$.
\end{theorem}

\begin{proof}
The claim is trivial for $k = 3$ since $H_c^{(3)}$ is the only fully symmetric nontrivial 3-peg Hanoi group. Assume $k \ge 4$, and let $G$ be a fully symmetric $k$-peg Hanoi group generated by $S$ that is not a subgroup of $H_c^{(k)}$. Then $S$ contains some $a \in S_k$ with two or more inactive pegs and nontrivial $\sigma_a$. We will show that $G$ is not contracting.

Let $p = |P_a|$. Relabel the pegs so that $Q_a = \{ p, p+1, \ldots, k-1\}$. Further relable pegs within $P_a$ so that $\sigma_a$ has the form $(0\ 1\ \cdots\ m)c_2 \ldots c_n$, a product of disjoint cycles. That is,
\[
a = (0\ 1\ \cdots\ m)c_2 \ldots c_n (1, 1, \cdots, 1, a, a, \cdots, a),
\]
where $a|_i = a$ for $i \ge p$. Take $b = (0\ p) \cdot a$, which is in $S$ by the symmetry assumption. Since each $c_i$, $2 \le i \le n$, is disjoint from $(0\ p)$,
\begin{align*}
b &= (p\ 1\ 2\ \ldots\ m)c_2 \ldots c_n (b, 1, \ldots, 1, 1, b, \ldots, b) \\
b^{-1} &= (m\ \ldots\ 2\ 1\ p)c_n^{-1} \ldots c_2^{-1} (b^{-1}, 1, \ldots, 1, 1, b^{-1}, \ldots, b^{-1}),
\end{align*}
where $b^{-1}|_i = b^{-1}$ for $i \in (0\ p)Q_a = \{0, p+1, \ldots, k-1\}$. Then $G$ contains
\[
ab^{-1} = (0\ 1\ p)(b^{-1}, a, 1, \ldots, 1, 1, ab^{-1}, \ldots, ab^{-1}),
\]
where $ab^{-1}|_i = ab^{-1}$ for $i \ge p+1$. Then $(ab^{-1})^3|_{p+1} = (ab^{-1})^3$ and $(ab^{-1})^3|_0 = ab^{-1}$, so Lemma \ref{noncontractioncondition} applies.
\end{proof}

Under weaker symmetry conditions, we obtain more contracting Hanoi groups. Identifying $\mathsf{X}_k$ with the vertices of a regular $k$-gon in the order 0 through $k-1$, we say that a Hanoi group is \emph{rotationally} or \emph{dihedrally symmetric} if its generating set is closed under the action of the corresponding symmetry group of $\mathsf{X}_k$. We do not have analogues of Theorem \ref{uniquecontracting} for these symmetries. We state one partial result.

\begin{definition} \label{Rkn}
Let $\underline{R}_{k, n} \subset S_k$, $2n + 1 \le k$, be the union of the identity automorphism and all automorphisms $a \in S_k$ for which
\begin{enumerate}
\item $Q_a \subset \{ m+1, m+2, \ldots, m+n \}$ for some $m$, and
\item whenever $i \in Q_a$, $\sigma_a$ also fixes every element of $\{ i - 1, i-2, \ldots, i-(n-1) \}$,
\end{enumerate}
where all computations are performed modulo $k$.
\end{definition}
In the cyclic arrangement, each peg has an increasing and a decreasing side. Condition 1 of Definition \ref{Rkn} says that the set of inactive pegs lie among $n$ adjacent pegs. Condition 2 says that $\sigma_a$ fixes the $n-1$ adjacent pegs on the decreasing side of each inactive peg.

\begin{proposition} \label{rotationalcontraction}
Every Hanoi group $(G, \mathsf{X}_k)$ generated by $S \subset \underline{R}_{k, n}$ is contracting.
\end{proposition}

\begin{proof}
We verify condition ($\ast$) from Theorem \ref{hanoicontraction2}. Take an arbitrary subset $T = \{ s_i \}_{i = 1}^M$ of $S$ with nonempty essential set $Q$, and consider any $j \notin \Fix(T)$. As usual, let $\sigma_i$ and $Q_i$ be the root permutation and the set of inactive pegs, respectively, of $s_i$. 

Without loss of generality, assume $n - 1 \in Q$. Then by the first condition on elements of $\underline{R}_{n,k}$, $Q_i \subset [0, 2n-2]$ for each $i$. Let $M_i$ be the largest element of $Q_i$, and choose $s_{min}$ so that $M_{min} \le M_i$ for all $i$. Note that $n-1 \le M_{min} \le M_i \le 2n-2$. Then $M_i + 1 \in Q_i$ and $M_i + 1 - (n - 1) \le n$, so by the second condition on elements of $\underline{R}_{n,k}$ each $\sigma_i$ fixes every element of $[n, M_{min}]$. Furthermore, since $n-1 \in Q_i$, each $\sigma_i$ fixes every element of $[0, n-1]$.

Since $j \notin \Fix(T)$, by the above we have $j \in [M_{min} + 1, k-1]$ and so  $\Orb_T(j) \subset [M_{min} + 1, k-1]$. Thus $\Orb_T(j) \cap Q_{min} = \emptyset$, as desired.
\end{proof}

Define $\overline{R}_{k, n}$ by replacing the set in the second condition of Definition \ref{Rkn} by $\{ i + 1, i + 2, i + (n - 1) \}$. The overline indicates that we now take the adjacent pegs on the increasing side. The analogue of Proposition \ref{rotationalcontraction} also holds with $\overline{R}_{k, n}$ instead of $\underline{R}_{k, n}$.

For example, Hanoi groups generated by the following sets are contracting by Proposition \ref{rotationalcontraction}:
\[
\{ 1 \} \cup \left\{
\begin{array}{cc}
a = (0\ 1) (1, 1, 1, a, a) & d = (3\ 4) (1, d, d, 1, 1) \\
b = (1\ 2) (b, 1, 1, 1, b) & e = (4\ 0) (1, 1, e, e, 1) \\
c = (2\ 3) (c, c, 1, 1, 1) &
\end{array}
\right\}
\]
\[
\{ 1 \} \cup \left\{
\begin{array}{cc}
a = (0\ 1) (1, 1, 1, a, a, 1) & d = (3\ 4) (d, d, 1, 1, 1, 1) \\
b = (1\ 2) (1, 1, 1, 1, b, b) & e = (4\ 5) (1, e, e, 1, 1, 1) \\
c = (2\ 3) (c, 1, 1, 1, 1, c) & f = (5\ 0) (1, 1, f, f, 1, 1)
\end{array}
\right\}.
\]
The first is a contracting 5-peg Hanoi group with rotational symmetry, while the second is a contracting 6-peg Hanoi group with dihedral symmetry. Moreover, unlike contracting $k$-peg Hanoi groups generated by subsets of $S_{k, k-2}$, $k > 3$ (recall Proposition \ref{HTnoncontracting}), these groups correspond to meaningful Hanoi games in the following sense.

\begin{proposition}
Let $G$ be a $k$-peg Hanoi group generated by $S$. The group generated by the set of all root permutations of the generators acts by permutation on $\mathsf{X}_k$. The Schreier graphs $\Gamma_n$, $n\ge 1$ are connected if and only if the action of $G$ on $\mathsf{X}_k$ is transitive.
\end{proposition}

\begin{proof}
Vertices $x, y \in \Gamma_n$ are connected if and only if $s_m \ldots s_2s_1(x) = y$ for some $s_i \in S \cup S^{-1}$. Since each $s_i$ has finite order, we can write this product without inverses. Thus $\Gamma_n$ is connected if and only if the legal moves of the associated Hanoi game allow us to change any legal $n$-disk configuration into any other legal $n$-disk configuration.

Suppose the action is not transitive, so that there exist $i$, $j \in \mathsf{X}_k$ such that no product $\sigma = \sigma_m \cdots \sigma_2\sigma_1$ of root permutations satisfies $\sigma(i) = j$. Then any configuration with the smallest disk on peg $j$ cannot be reached from any configuration with the smallest disk on peg $i$.

Now assume the action is transitive. The case $n = 1$ is immediate. Assume the claim for $n - 1$, and consider any two legal $n$-disk configurations. Suppose the largest disk needs to move from peg $i$ to peg $j$. By transitivity, $\sigma_m \ldots \sigma_2\sigma_1(i) = j$ for some product of root permutations. Let $s_i$ be some generator with root permutation $\sigma_i$. Repeat the following for $1 \le i \le m$: \emph{Use the $n-1$ case to move the $n-1$ smaller disks onto an inactive peg of $s_i$. Then use $s_i$ to move the largest disk.} After the $m$-th step, the largest disk will be on peg $j$. We then use the $n-1$ case to move the smaller disks to their desired pegs.
\end{proof}

We end with several open questions and thoughts.

\begin{itemize}
\item Theorem \ref{uniquecontracting} holds with the alternating group on $\mathsf{X}_k$; since $(0\ p\ p+1) \cdot a = (0\ p) \cdot a$, the proof only requires 3-cycles. Is there a weaker or more natural set of constraints under which $H_c^{(k)}$ remains the unique maximal contracting group?
\item There are non-Hanoi automorphisms that may still be interpreted in terms of the Hanoi Towers game. For example, consider
\[
a = (0\ 1\ 2) (1, 1, 1, b) \quad b = (0\ 2\ 1) (1, 1, 1, a).
\]
Then $ab = (1, 1, 1, ba)$ and $ba = (1, 1, 1, ab)$, so $ab = 1$. In terms of the game, both $a$ and $b$ start from the smallest disk and keep track of parity the number of disks on peg 3 until it encounters a disk on pegs 0, 1, or 2. If the parity is odd, $a$ applies $(0\ 1\ 2)$; if even, $(0\ 2\ 1)$. $b$ does the same with the permutations exchanged. How can we characterize all automorphisms that correspond to legal moves of the game in some reasonably defined sense? Which groups with such generators are contracting?
\item Can we apply similar symmetry conditions to other self-similar groups to obtain criteria on their contraction or non-contraction?
\end{itemize}

\appendix
\section{Appendix}\label{Appendix}
This paper began by examining existing variations of the Hanoi Towers game. Two such modifications were the Hanoi Networks, HN3 and HN4, which were introduced in the physics literature \cite{BGG08b} as examples of regular networks with small-world properties. For more on networks with this small-world property, see \cite{BGG08}, \cite{NW99}, and \cite{WS98}. We will show that while HN3 can be related peripherally to the automaton representing the Hanoi Towers game, these networks have no direct connection to the original Hanoi Towers game.

\subsection{Construction of HN3/HN4}
The construction of both networks is only partially inspired by the traditional 3-peg Hanoi Tower game. 

\begin{definition}
Define \emph{$S_n$} as the sequence of disks moved in the optimal solution of the $n$-disk game, where disks numbered $1$ through $n$ from smallest to largest. That is if you move disk1 on the first move, the first element in the sequence is 1.
\end{definition}

For example, the optimal solution for the 1-disk game simply moves disk 1 to the desired peg, so that
\[
S_1 = ( 1 ).
\]
For two disks, we first move disk 1 to the third peg, disk 2 to the goal peg, then disk 1 above disk 2:
\[
S_2 = ( 1, 2, 1 ).
\]

It is well known (see, for example, \cite{BH99}) that the optimal solution is given recursively; to move $n$ disks from peg 1 to peg 2, we first move the $n-1$ smaller disks to peg 3, move disk $n$ to peg 2, and finally move the $n-1$ smaller disks from peg 3 to peg 2. This means that $S_{n+1}$ is obtained by concatenating two copies of $S_n$, with a $n+1$ term between them. We can therefore define the following:

\begin{definition}
Define $S$ as the unique infinite sequence having each $S_n$ as a prefix:
\[
S = ( 1, 2, 1, 3, 1, 2, 1, 4, 1, 2, 1, 3, 1, 2, 1, 5, \ldots )
\]
\end{definition}

A \emph{network} is a graph with a length assigned to each edge. Vertices of a network are called \emph{nodes}. We use $S$ to construct two networks that incorporate this model of the Hanoi Towers game and also have small world properties.

\begin{figure}[htb]
\begin{center}
\includegraphics[width=80mm]{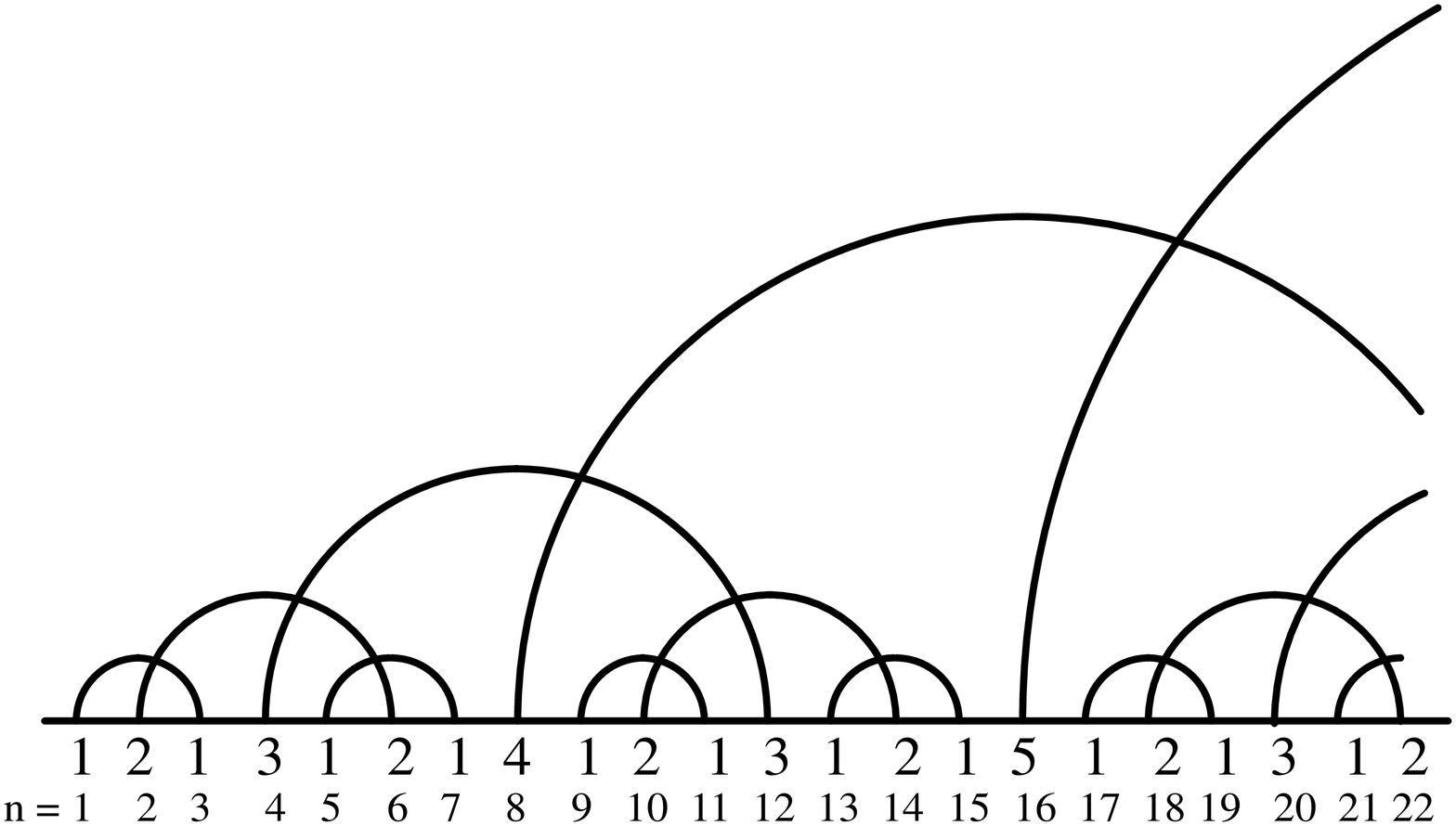} \vskip 10mm
\includegraphics[width=80mm]{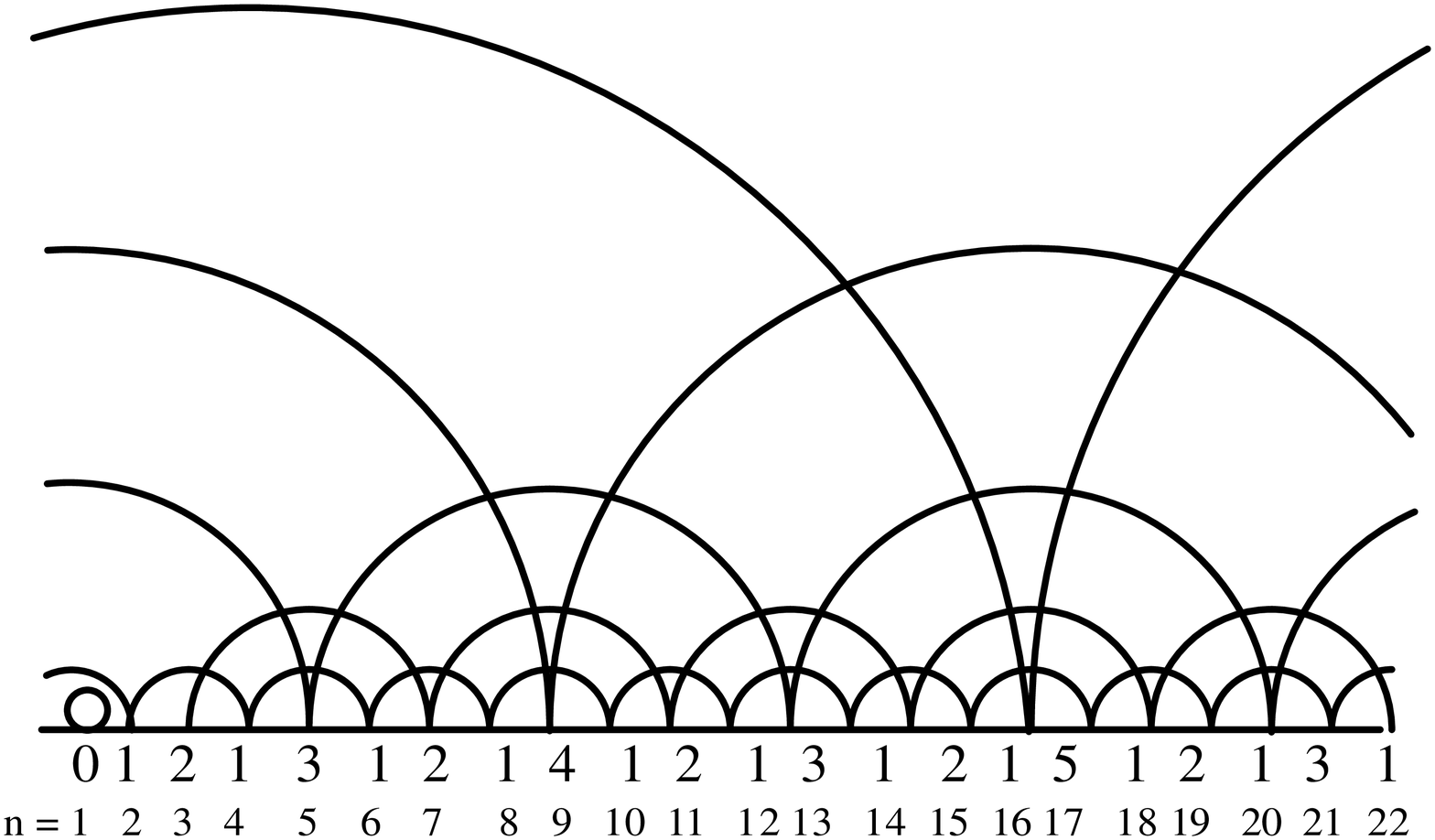}
\end{center}
\caption{The networks HN3 and HN4}\label{HN3andHN4}
\end{figure}

\begin{definition}
The HN3 network has nodes that are identified with the positive integers. Nodes corresponding to adjacent integers are connected; these edges form the \emph{backbone} of this network. Label node $n$ by the $n$-th element of $S$, called the \emph{disk number} of the node. Nodes with the same disk number $i$ are connected if and only if there exists a node labeled $i+1$ in between and there does not exist a node labeled $j$, where $j \ge {i+2}$, in between. We call these connections \emph{long-distance jumps}. An edge connecting nodes $n$ and $m$ is given length $|n - m|$. See Figure \ref{HN3andHN4} for consistency with Definition \ref{def:hn4}.
\end{definition}

Each node is thus connected at least to its adjacent nodes on the backbone by edges of length 1. In addition, certain nodes with disk number $i$ are connected by long-distance jumps of length $2^i$. For example, nodes 1 and 3, the first two nodes with disk number 1, are connected, as are 5 and 7, 9 and 11, etc; among nodes with disk number 2, nodes 2 and 6 are connected, then 10 and 14, etc.; and so on for every disk number.

\begin{definition}\label{def:hn4}
The network HN4 is the network HN3 with additional connections made between a node labeled $i$ and the first nodes labeled $i$ to its left and to its right. In addition, the backbone is extended to all integers, and the network is made symmetric over 0. The special point 0 is connected to itself by a loop. See Figure \ref{HN3andHN4}.
\end{definition}

For example, for disk number 1, HN4 retains all the long-distance jumps in HN3, and in addition connects $n=3$ to $n=5$, $n=7$ to $n=9$, etc.  to create a 4-regular network. In both HN3 and HN4, the metric is given by the distance along the backbone (not using any long-distance jumps). 

\subsubsection{Hanoi Tower Automata} \label{HanoiTowerAutomaton}
The game also gives rise to a series of finite networks corresponding to the finite state automata of the $n$-disk game. 

\begin{definition}
The network \emph{$H_n$} is constructed to represent the $n$-disk game so that each vertex of the automaton is associated to an $n$-letter word $x_1 x_2 x_3 \ldots\ x_n$, $x_i \in \{0, 1, 2\}$. This word indicates the peg number of each disk, starting from the smallest disk. Two vertices are connected by an edge if a legal move allows the player to move between the corresponding states. Every edge is given the same length, chosen so that the shortest path between $0^{n}$ and $1^{n}$ has length 1.
\end{definition}

For notational simplicity, we write $i^n$ for the $n$-letter word $i \ldots ii$.

Note that the vertices of $H_n$ correspond to the possible states of the game, i.e. the configuration of $n$ disks. While each $H_n$ has a natural metric independent of the ambient space, it is useful to have an embedding of these networks, both as a visualization and in order to obtain some limiting object. We find that a particular recursive construction allows us to embed each $H_n$ in the plane so that every edge has the same length.

\begin{figure}[htbp]
\begin{center}
\includegraphics[width=60mm]{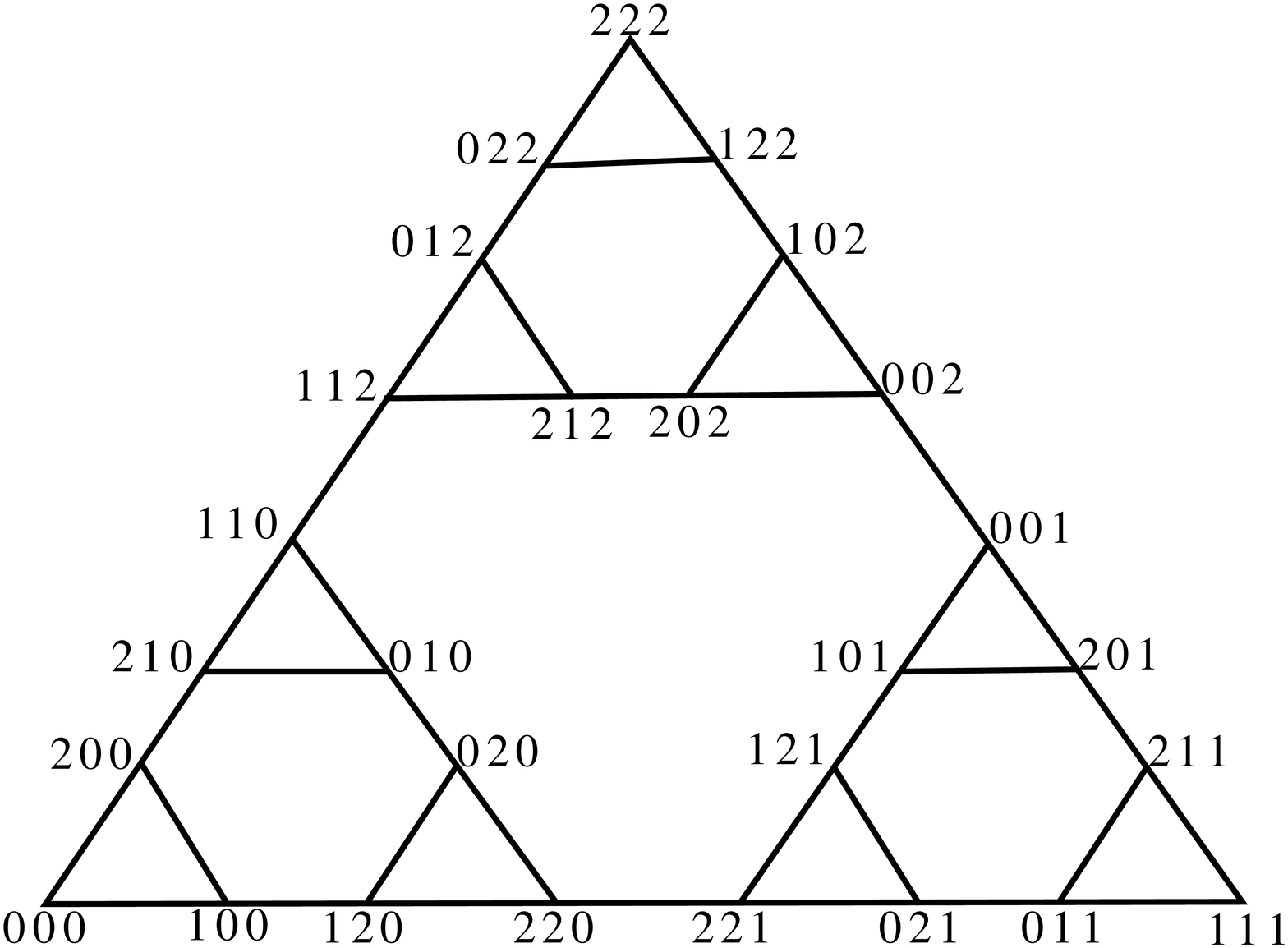}
\end{center}
\caption{H$_3$}
\label{LabeledH_3}
\end{figure}

The network $H_1$ has three states: 0, 1, and 2. We arrange these in an equal triangle; to fix ideas, place 0 at the bottom left, 1 at the bottom right, and 2 at the top. Every $H_n$ will contain a triangle with vertices $0^n$, $1^n$, and $2^n$ in this same orientation. For the $n$-disk game represented by $H_n$, $n \ge 2$, ignoring the largest disk results in the Hanoi Tower game for the $n-1$ smaller disks. 
As a result, $H_n$ contains three transformed copies of $H_{n-1}$, each representing the $(n-1)$-disk game with the largest disk on one of the three pegs. The three copies are transformed in the following way. For the copy of $H_{n-1}$ corresponding to the largest disk being on peg $i$, where $i=0,1,2$, reflect $H_{n-1}$ with its labels across the line through $i^{n-1}$ and the pidpoint of the other two vertices, and append $i$ to each label.
 Now, translate these three copies apart so that the new labels $0^n$, $1^n$, and $2^n$ form the vertices of a larger triangle. The only missing edges of $H_n$ are those that correspond to moving the largest disk. This is only possible when the $n-1$ smaller disks are on the same peg, different from the one with the largest disk. It easily follows that there are three more edges: between $1^{n-1}2$ and $1^{n-1}0$; $2^{n-1}0$ and $2^{n-1}1$; $0^{n-1}1$ and $0^{n-1}2$. With appropriate translations of the smaller networks, these three edges connect the three triangles to complete the outer edge of the larger triangle. Moreover, their lengths can be chosen to be the same as the edges in each copy of $H_{n-1}$, giving the desired embedding of $H_n$. Compare this with Figure \ref{LabeledH_3}.

\subsection{Relation between HN3 and Hanoi Tower Automaton}
In the definitions of HN3 and HN4 the long-distance jumps were added to create a small-world property and as such have no direct relation to the Hanoi Tower game, there is no reason to expect a connection between these and $H_n$. However, the recursive nature of HN3 allows us to obtain it as a subnetwork of the automaton network.

\begin{figure}[htbp]
\begin{center}
\includegraphics[width=100mm]{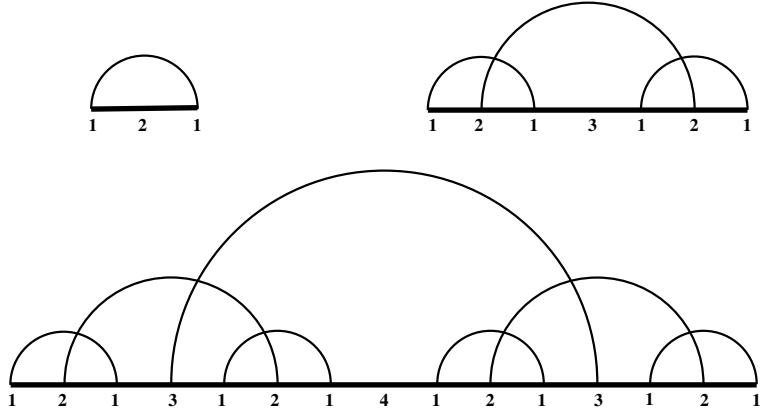}
\end{center}
\caption{The first three levels of HN3$_n$}
\label{LabeledHN3(level1-3)}
\end{figure}

\begin{figure}[htbp]
\begin{center}
\includegraphics[width=45mm]{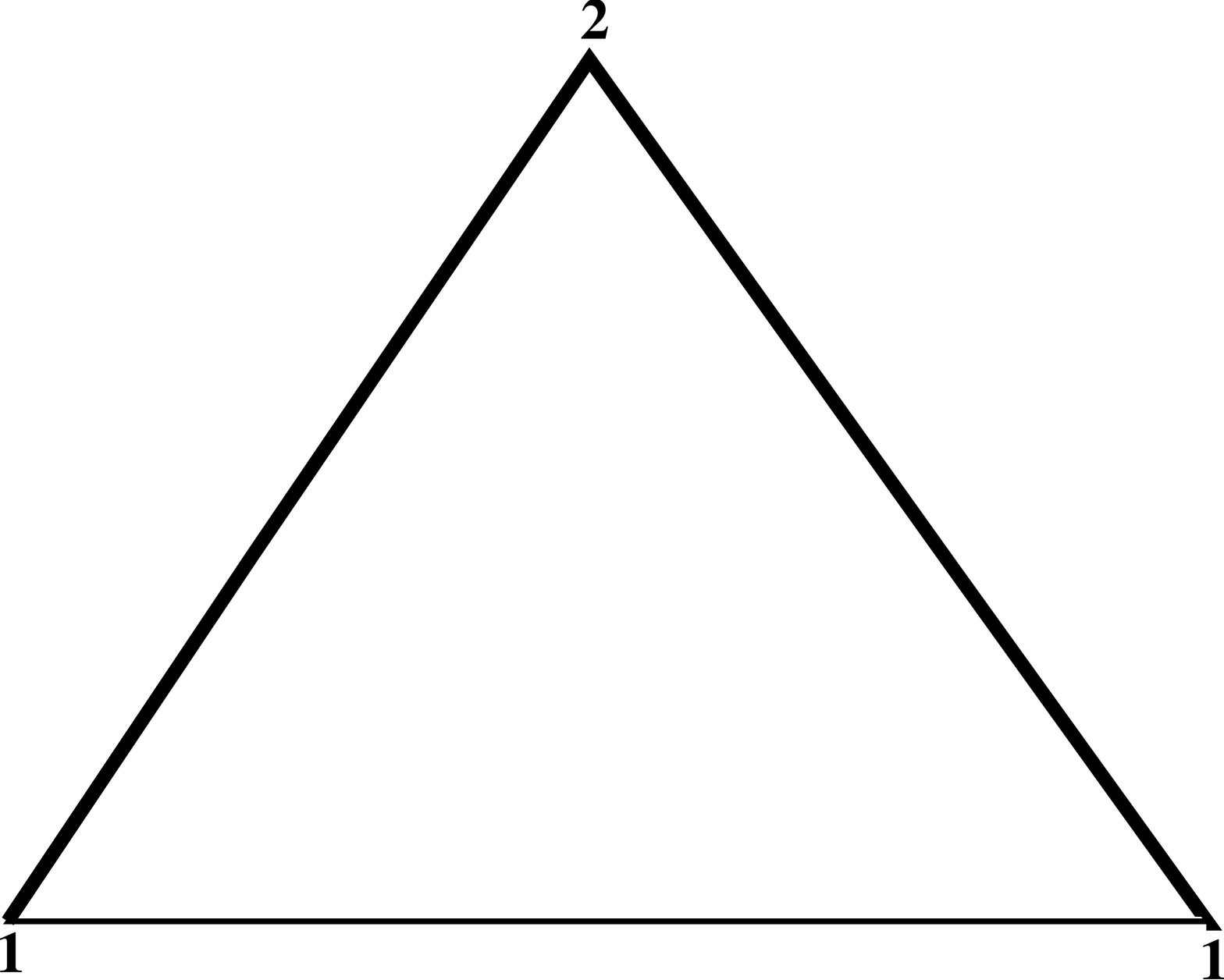} \hskip 5 mm
\includegraphics[width=45mm]{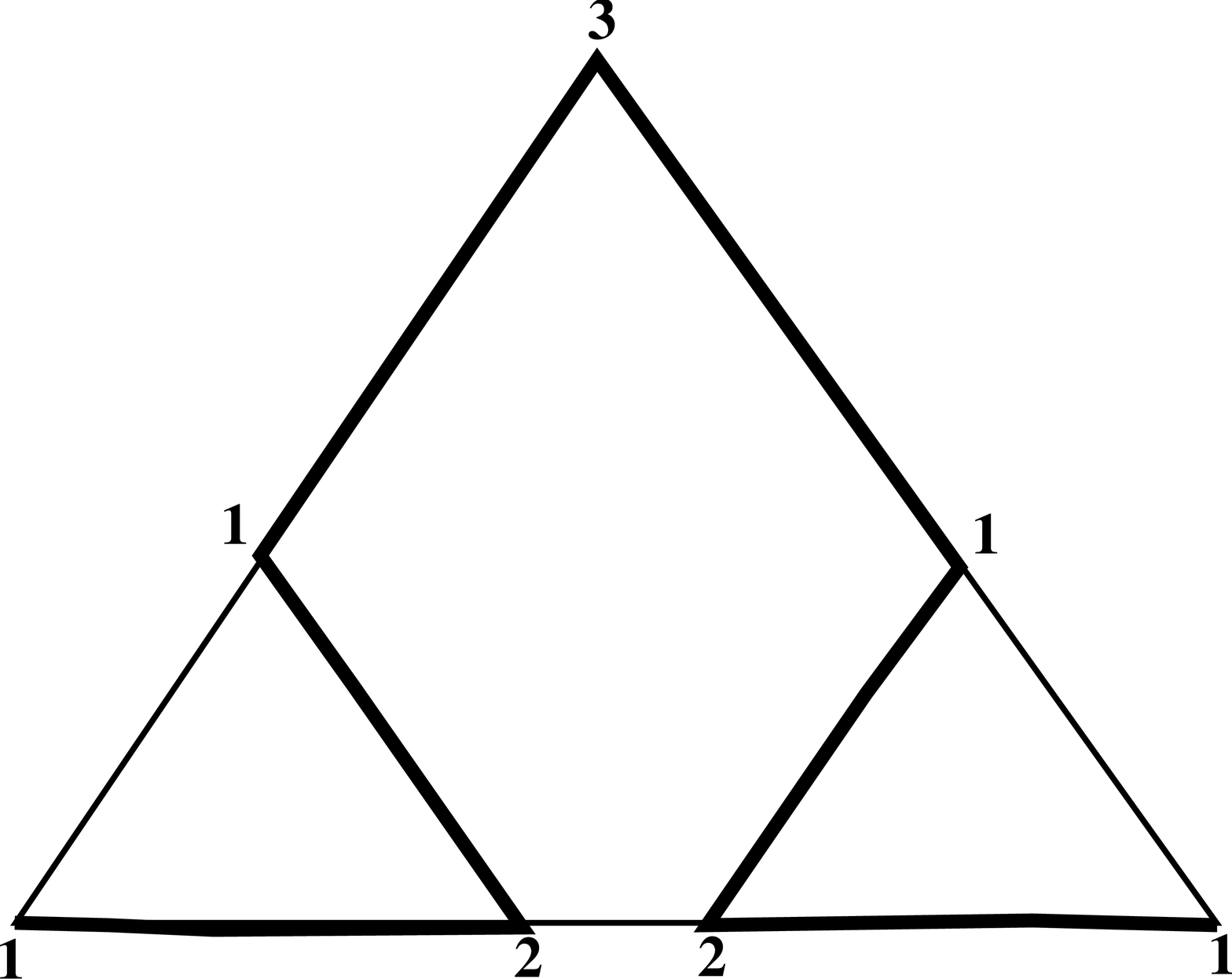} \vskip 5 mm

\includegraphics[width=50mm]{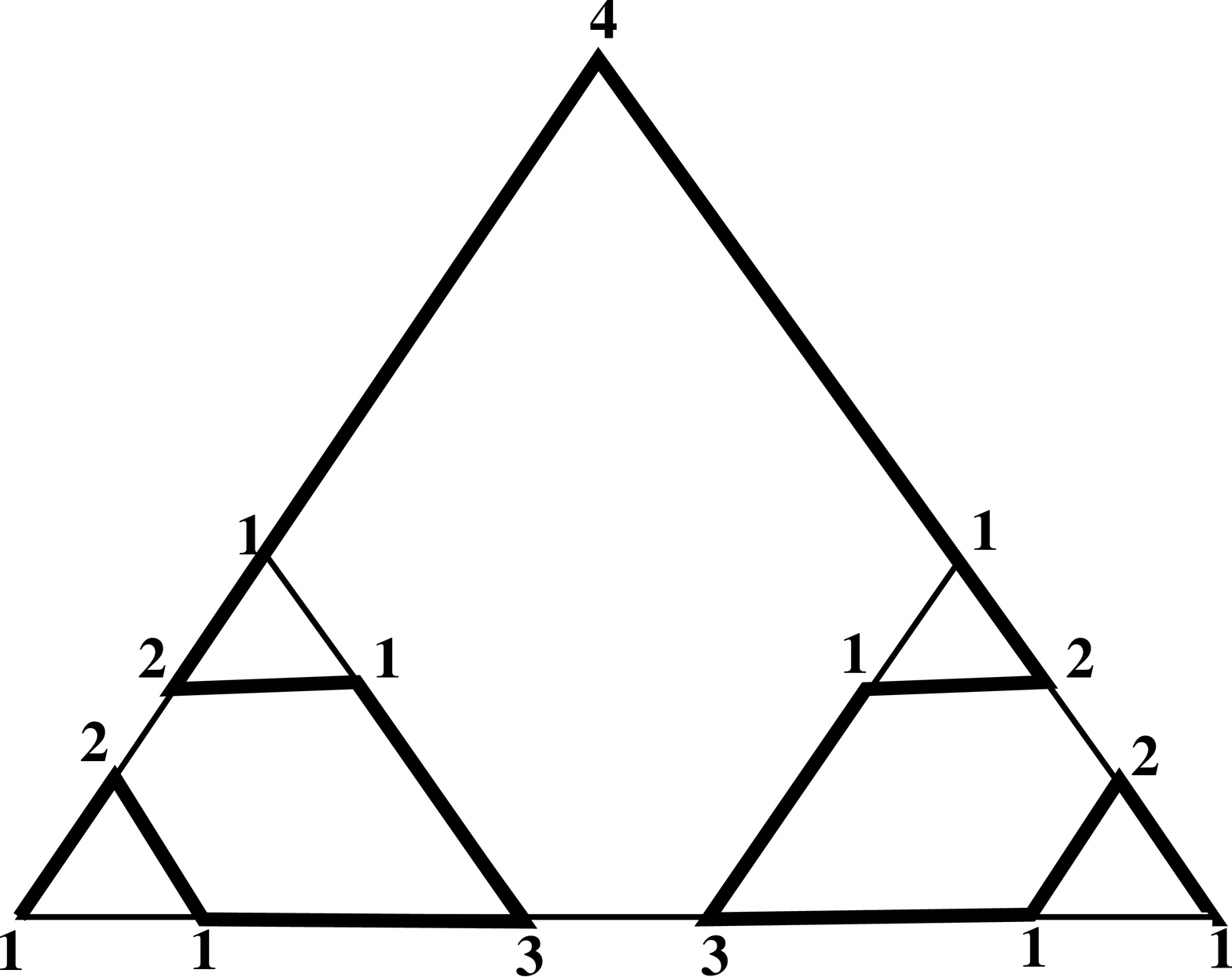}
\end{center}
\caption{The first three level of $H'_n$. The darkened lines correspond to the backbone of HN3$_n$ while the lighter lines correspond to the long-distance jumps.} \label{H'_n(level1-3)}
\end{figure}

\begin{theorem}\label{subnetworkthm}
Let HN3$_n$ be the subnetwork of HN3 consisting of the first $|S_n|$ vertices and edges connecting them. Then HN3$_n$ can be obtained as a graph minor $H'_n$ of $H_n$.
\end{theorem}
\begin{proof}
For $n=1$, set $H_1' = H_1$. Then HN3$_1$ is clearly isomorphic to $H'_1$ as graphs. We develop the isomorphism between HN3$_n$ and $H'_n$ inductively so that the middle node of HN3$_n$, with disk number $n + 1$, corresponds to $2^n$ in $H_n$, and the end nodes with disk number $1$ correspond to $0^n$ and $1^n$.

Fix $n \ge 2$. Analogous to the way $S_n$ is obtained from two copies of $S_{n-1}$, HN3$_n$ can be obtained from two copies of HN3$_{n-1}$, an extra node of disk number $n$, and three new edges: two edges to join the backbone of each copy of HN3$_{n-1}$ to the new node, and a long-distance jump connecting the two nodes of disk number $n-1$ at the center of each copy of HN3$_{n-1}$. To construct $H'_n$, we proceed as in the recursive construction of $H_n$, but create only two transformed copies of $H'_{n-1}$, say for $i = 0$ and $1$, using the reflections and label appending explained in \ref{HanoiTowerAutomaton}. These correspond to the two copies of HN3$_{n-1}$. In place of the third copy of $H'_{n-1}$, $H'_n$ has the single node $2^n$, corresponding to the added middle node in HN3$_n$ of disk number $n$. The nodes $1^{n-1}0$ and $0^{n-1}1$ are connected to the node $2^n$, corresponding to the two new edges that complete the backbone of HN3$_{n-1}$. The nodes $2^{n-1}0$ and $2^{n-1}1$ are also connected, corresponding to the added long-distance jump in HN3$_n$.

These correspondences show that the isomorphism of HN3$_{n-1}$ and $H'_{n-1}$ implies the isomorphism of HN3$_n$ and $H'_n$, so by induction, HN3$_n = H'_n$ for all $n \ge 1$. By construction, $H'_n$ can be regarded as a subgraph of $H_n$ in a natural way. Figures \ref{LabeledHN3(level1-3)} and \ref{H'_n(level1-3)} show the first three levels of $HN3_n$ and $H'_n$, respectively, with the backbone of $HN3_n$ and the corresponding edges in $H'_n$ thickened.
\end{proof}

In this correspondence, every long-distance jump in HN3 receives one edge-length in any $H_n'$. All edges of the backbone of HN3 connect a node of disk number 1 to another node, say of disk number $k$. These edges receive $k-1$ edge-lengths in any $H_n'$. As a result, the ratio of distortion of the metric is unbounded, but only as disk number also goes to infinity. Since higher disk numbers appear further away from the origin in HN3, the distortion is controlled for any finite section of the networks. However, because of difference in metric, it is impossible to extend the correspondence to the limiting space. More precisely, we have the following theorem:

\begin{theorem}
There exists a Lipschitz embedding of HN3$_n$ into $H_n$ that is not uniformly bi-Lipschitz in $n$.
\end{theorem}
\begin{proof}
The subnetwork $H'_n$ is an embedding of HN3$_n$ into $H_n$ that is Lipschitz since all edges in $H_n$ have length less than 1 and all distances in HN3$_n$ are greater than or equal to 1. However, the lower Lipschitz bound is smaller than $2^{-n}$. Thus the embedding is bi-Lipschitz for each $n$, but not uniformly.
\end{proof}

The network $H'_n$ has an interpretation as an automaton for its own game. Start with the same labeling 0, 1, and 2 for $H'_1$. At each iteration, follow the usual labeling rule for the two copies of $H'_{n-1}$, but give the label $2$ to the extra node. Then each label in $H'_{n}$ is, read from right to left, a binary number which ends after $n$ digits or at the first appearance of 2. In terms of the Hanoi Towers game, this corresponds to ignoring disks smaller than the largest disk at peg 2, if there are any. All states in the Hanoi Towers game that agree from the largest disk to the largest disk at peg 2 represent the same state in $H'_n$, just as we collapsed the copy of $H'_{n-1}$ with the largest disk on peg 2 to a single vertex in the recursive construction. It also follows that the disk number derived from HN3$_n$ has a natural interpretation: a node in $H'_n$ with disk number 1 represents a single state in $H_n$, and a node with disk number $k > 1$ represents the result of collapsing $3^{k-2}$ states in $H_n$.


\small
\begin{thebibliography}{99}

\bibitem
{physics}  E. Akkermans, G. Dunne, \& A. Teplyaev, \emph{Physical Consequences of Complex Dimensions of Fractals}  (2009).

\bibitem
{Betal08}N. Bajorin, T. Chen, A. Dagan, C. Emmons, M. Hussein, M. Khalil, P. Mody, B. Steinhurst, and A. Teplyaev, \emph{Vibration modes of 3$n$-gaskets and other fractals}, J. Phys. A \textbf{41} (2008), no. 1, 015101, 21 pp.

\bibitem
{BGG08}S. Boettcher, B. Gon\c{c}alves and H. Guclu, \emph{Hierarchical Regular Small-World Networks}, J. Phys. A \textbf{41} (2008), no. 25, 252001, 7 pp.

\bibitem
{BGG08b}Stefan Boettcher and Bruno Gon\c{c}alves, \emph{Anomalous Diffusion on the Hanoi Networks}, arXiv:0802.2757v2.

\bibitem
{BN03}Evgen Bondarenko and Volodymyr Nekrashevych, \emph{Post-critically finite self-similar groups}, Algebra Discrete Math. \textbf{2} (2003), no. 4, 21--32.

\bibitem
{BH99}J.-P. Bode and A. M. Hinz, \emph{Results and open problems on the Tower of Hanoi}, Congr. Numer. \textbf{139} (1999), 113-–122.

\bibitem
{FS92}M. Fukushima and T. Shima, \emph{On a spectral analysis for the Sierpi\'{n}ski gasket}, Potential Anal. \textbf{1} (1992), no. 1, 1--35.

\bibitem
{Gri80}R. I. Grigorchuk, \emph{On Burnside's problem on periodic groups (Russian)}, Funktsional. Anal. i Prilozhen. \textbf{14} (1980), no. 1, 53--54. 

\bibitem
{GS06}Rostislav Grigorchuk and Zoran \v{S}un\'{i}k, \emph{Asymptotic aspects of Schreier graphs and Hanoi Towers groups}, C. R. Math. Acad. Sci. Paris \textbf{342} (2006), no. 8, 545–-550.

\bibitem
{GS08}Rostislav Grigorchuk and Zoran \v{S}un\'{i}k, \emph{Schreier spectrum of the Hanoi Towers group on three pegs}, \emph{Analysis on graphs and its applications}, 183--198, Proc. Sympos. Pure Math, vol. 77, Amer. Math. Soc., Providence, RI, 2008.

\bibitem{kajino} N. Kajino, \emph{Cell-counting dimension and spectral asymptotics for Laplacians on self-similar sets} (2009).

\bibitem
{Kig01}Jun Kigami, \emph{Analysis on Fractals}, Cambridge Tracts in Mathematics, vol. 143. Cambridge University Press, Cambridge, 2001.

\bibitem
{Kig95}Jun Kigami, \emph{Hausdorff dimensions of self-similar sets and shortest path metrics}, J. Math. Soc. Japan \textbf{47} (1995), no. 3, 381--404.

\bibitem
{KL93}Jun Kigami and Michel L. Lapidus, \emph{Weyl's problem for the spectral distribution of Laplacians on p.c.f. self-similar fractals}, Comm. Math. Phys. \textbf{158} (1993), no. 1, 93--125.

\bibitem
{Mor46}P. A. P. Moran, \emph{Additive functions of intervals and Hausdorff measure}, Proc. Cambridge Philos. Soc. \textbf{42} (1946), 15--23.

\bibitem
{Nek05}Volodymyr Nekrashevych, \emph{Self-similar groups}, Mathematical Surveys and Monographs, 117, Amer. Math. Soc., Providence, RI, 2005.

\bibitem
{NT08}Volodymyr Nekrashevych and Alexander Teplyaev, \emph{Groups and analysis on fractals}, \emph{Analysis on graphs and its applications}, 143--180, 
Proc. Sympos. Pure Math., vol. 77, Amer. Math. Soc., Providence, RI, 2008.

\bibitem
{NW99}M. E. J. Newman and D. J. Watts, \emph{Renormalization group analysis of the small-world network model}, Physics Letters A \textbf{263} (1999), 341-346.

\bibitem
{Str06}Robert S. Strichartz, \emph{Differential Equations on Fractals. A Tutorial}, Princeton University Press, Princeton, NJ, 2006.

\bibitem
{WS98}D. J. Watts and S. H. Strogatz, \emph{Collective dynamics of `small-world' networks}, Nature \textbf{393} (1998), 440--442.

\end{thebibliography}
\end{document}